\documentclass[12pt]{amsart}
\usepackage{amssymb}
\newtheorem{thm}{Theorem}[section]
\newtheorem{crl}[thm]{Corollary}

\newtheorem{lm}[thm]{Lemma}

\newcommand{\der }{\partial}

\newcommand{\M}{\mathcal M}
\newcommand{\D}{\mathcal D}

\begin{document}

\title
{$N$-commutators of vector fields}
\author{A.S. Dzhumadil'daev}
\address{Institut des Hautes \'Etudes Scientifiques, Bures-sur-Yvette
\newline
Institute of Mathematics, Academy of Sciences of Kazakhstan,
{\hbox{Almaty}}} \email{askar@math.kz}

\maketitle

\begin{abstract}
The question whether the $N$-commutator
$s_N(X_1,\ldots,X_N)= \sum_{\sigma\in
Sym_N}sign\,\sigma\,X_{\sigma(1)}\cdots X_{\sigma(N)}$ is a well
defined operation on $Vect(n)$ is studied.
It is if $N=n^2+2n-2.$ A theory of
$N$-commutators with emphasis on $5$-and $6$-commutators on
two-dimensional manifolds is developed.
 Divergenceless vector fields $Vect_0(2)$ under $2$-~commutator and
$5$-commutator generate a $sh\mbox{-Lie}$ algebraic structure in
the sense of Stasheff. Similar results hold for $Vect(2)$ with $2$- and
$6$-commutators.
\end{abstract}

\tableofcontents
\section{Introduction}
Let $M=M^n$ be a $n$-dimensional smooth flat manifold, $U$ the algebra of
smooth functions on $M$  and $Dif\!f(n)$ an algebra of
$n$-dimensional differential operators $\sum_{\alpha\in {\bf
Z}^n}\lambda_{\alpha}u_{\alpha}\der^{\alpha}$ where
$\der^{\alpha}=\prod_{i=1}^n\der_i^{\alpha_i},$
$\der_i=\der/\der\,x_{i},$ $\alpha=(\alpha_1,\ldots,\alpha_n),$
$\alpha_i\in {\bf Z},$ and
$u_{\alpha}\in U.$ The subspace of differential operators of first order
$Vect(n)=\{\sum_{i=1}^nu_i\der_i\}$
can be interpreted as a space of vector fields on $M$.

The algebra  $Dif\!f(n)$ is associative but its subspace $Vect(n)$ is
not closed under composition.
If $X=u_i\der_i, Y=v_j\der_j,$ then
$$X\cdot Y= u_i\der_i(v_j)\der_j+u_iv_j\der_i\der_j,$$
has a quadratic differential part $u_iv_j\der_i\der_j.$
But $Vect(n)$ is  closed under commutator. We see that
$$Y\cdot X=v_j\der_j(u_i)\der_i+v_ju_i\der_j\der_i,$$
has also a quadratic differential part that is cancelled by the
quadratic differential part of $X\cdot Y.$ Since $M$ is flat,
$\der_i\der_j=\der_j\der_i$ and
$$[X,Y]=X\cdot Y-Y\cdot X=u_i\der_i(v_j)\der_j
-v_j\der_i(u_i)\der_i\in Vect(n).$$

This is a well known fact, that in terms of skew-symmetric polynomials means
the following.
Let $Sym_k$ be a permutation group,
$sign\,\sigma$ the sign of a permutation $\sigma\in Sym_k$ and
$$s_k(t_1,\ldots,t_k)=\sum_{\sigma\in Sym_k}sign\,\sigma\,t_{\sigma(1)}\cdots t_{\sigma(k)}$$ standard skew-symmetric associative polynomial.
Then we can substitute instead of $t_i$ any differential operator
from $Dif\!f(n).$ Closedness of $Vect(M)$ means that $s_2$ is a well
defined operator, i.e.,
$$s_2(X,Y)\in Vect(M),$$
for any $X,Y\in Vect(M).$

Now we consider the question: does there exist $k>2,$ such that
$s_k$ is also a well defined operation on $Vect(n)$? In other
words, does there exist some $k>2,$ such that
$$s_k(X_1,\ldots,X_k)\in Vect(n),$$
for any $X_1,\ldots,X_k\in Vect(n)$? Since $X_i=\sum_{j=1}^nu_{ij}\der_j,
i=1,2,\ldots,k,$ are differential operators of first order,
$X_{\sigma(1)}\cdots X_{\sigma(k)}$ are in general differential operators
of order $k$ and $s_k(X_1,\ldots,X_k)$ is
a differential operator of order $k.$

Surprisingly, for some special $k=k(n)$ it might happen that all
higher degree differential parts of $s_k(X_1,\cdots,X_k),$ like
quadratic differential part of $s_2,$ can be  cancelled for any
$X_1,\ldots,X_k\in Vect(n),$ but the first order part not.

Let us  demonstrate this observation on the
space of vector fields on 2-dimensional manifolds. We prove that
for $n=2$ $s_6$ is a well defined non-trivial $6$-linear map on
$Vect(2):$
$$\forall X_1,\ldots,X_6\in Vect(2)\Rightarrow s_6(X_1,\ldots,X_6)\in Vect(2),
$$
and
$$s_6(X_1,\ldots,X_6)\ne 0, \mbox{for some}\quad X_1,\ldots,X_6\in
Vect(2). $$
The number $6$ here can not be improved:
$$k>6\Rightarrow s_k(X_1,\ldots,X_k)=0, \quad
\forall X_1,\ldots,X_k\in Vect(2),$$
and $s_k(X_1,\ldots,X_k), 2<k<6,$ has a non-trivial quadratic differential
part for some $X_1,\ldots,X_k\in Vect(2).$

We also prove that, for the subspace of divergenceless vector fields
$Vect_0(2),$ the number $6$ can be improved. It has a non-trivial
$5$-commutator:
$$\forall X_1,\ldots,X_5\in Vect_0(2)\Rightarrow s_5(X_1,\ldots,X_5)\in
Vect_0(2),$$
and
$$ s_5(X_1,\ldots,X_5)\ne 0, \mbox{for some}\;\;X_1,\ldots,X_5\in Vect_0(2). $$
Moreover, the number $5$ here can not be changed:
$$k>5\Rightarrow s_k(X_1,\ldots,X_k)=0, \forall X_1,\ldots,X_k\in Vect_0(2).$$
If $2<k<5,$ then $s_k(X_1,\ldots,X_k)$ has a non-trivial quadratic
differential part for some $X_1,\ldots,X_k\in Vect_0(2).$

So, the vector space $Vect(2)$ can be endowed with a structure of
algebra under 2-commutator $s_2,$ usually denoted by $[\;,\;],$
and $6$-commutator $s_6.$ Similarly,
$Vect_0(2)$ can be endowed by a structure of Lie algebra under
$2$-commutator $s_2$ and $5$-commutator $s_5.$  These commutators have some
nice properties.

Let $A=Vect(2)$ or $Vect_0(2)$ and $X,Y, X_1,X_2,\ldots\in A.$
It is well known that an adjoint derivation $ad\,X,$ defined by
$(Y)ad\,X=[Y,X],$ is a derivation of $A.$ It is also well known that vector
fields can be interpreted as differential operators of first
order and in this sense $Vect(n)$ is a subspace of a space of all
differential operators $Dif\!f(n).$ The space of differential
operators have a composition operation, denoted by $\cdot$, and
under this operation $(Dif\!f(n),\cdot)$ is an associative algebra.
The $2$-commutator $[X,Y]=X\cdot Y-Y\cdot X$ as we mentioned
can be represented in terms of right-symmetric multiplication
$X\cdot Y-Y\cdot X= X\circ Y-Y \circ X.$ As it turns out, these
facts have analogues for $5$- and $6$-commutators.

The following Leibniz rule holds:
$$[X,s_5(X_1,\ldots,X_5)]=\sum_{i=1}^5s_5(X_1,\ldots,X_{i-1},[X,X_i],X_{i+1},\ldots,X_5),$$
for any $X,X_1,\ldots,X_5\in Vect_0(2).$
To calculate the $5$-commutator of $X_1,\ldots,X_5,$
one can use right-symmetric multiplication:
$$\sum_{\sigma\in Sym_5}sign\,\sigma\,
X_{\sigma(1)}\cdot X_{\sigma(2)}\cdot X_{\sigma(3)}\cdot X_{\sigma(4)}\cdot X_{\sigma(5)}=$$
$$\sum_{\sigma\in Sym_5}sign\,\sigma\,
(((X_{\sigma(1)}\circ X_{\sigma(2)})\circ X_{\sigma(3)})\circ
X_{\sigma(4)})\circ X_{\sigma(5)}.$$
In other words,
$$s_5(X_1,\ldots,X_5)=s_5^{rsym.r}(X_1,\ldots,X_5), \quad
\forall X_1,\ldots,X_5\in Vect_0(2).$$ Using right-symmetric
multiplication makes the calculation of $5$-commutators  easier
than the calculation of a composition of $5$ differential operators.
The $5$-commutator satisfies the following condition
$$\sum_{\sigma\in Sym_5}sign\,\sigma\,s_5(X_{\sigma(1)},X_{\sigma(2)},X_{\sigma(3)},X_{\sigma(4)},
s_5(X_{\sigma(6)},X_{\sigma(7)},X_{\sigma(8)},X_{\sigma(9)},X_0))=0,$$
for any $X_0,X_1,\ldots,X_9\in Vect_0(2).$ We call this identity the
{\it $4$-left commutativity  identity.} The algebra $Vect_0(2)$
endowed with multiplications $\omega_1,\omega_2, \omega_3,\ldots,$
such that $\omega_i=s_i$ for $i=2,5$ and $\omega_i=0$ for $i\ne
2,5,$ is a {\it $sh$-algebra} or {\it strongly homotopical algebra}
in the sense of Stasheff \cite{StashefLada}. In particular, $s_5$ is
a $5$-cocycle of the adjoint module of the Lie algebra
$(Vect_0(2),[\;,\;]).$ Moreover, the following
relations hold:
$$\sum_{1\le i<j\le 6}(-1)^{i+j}s_5^{rsym.r}([X_i,X_j],X_1,\ldots,\hat{X_i},\ldots,\hat{X_j},\ldots,X_6)=0,$$
$$\sum_{i=1}^6(-1)^i[X_i,s_5^{rsym.r}(X_1,\ldots,\hat{X_i},\ldots,X_6)]=0,$$
for any $X_1,\ldots,X_6\in Vect_0(2).$

For any two vector fields their commutator is once again a vector
field. One can repeat  this procedure $k-1$ times and construct from
any $k$ vector fields a new vector field. This can be done in many
ways. One can get a linear combination of such commutators. So, in
general there are many ways to construct invariant $k$-operation
on $Vect_0(n)$ or $Vect(n).$ We consider invariant operations
obtained in such ways as standard. Call a $k$-invariant
non-standard operation on vector fields a {\it primitive
$k$-commutator}. We prove that  the $5$-commutator and the
$6$-commutator are primitive, i.e., they can not be obtained from
commutators by such a procedure.

Any divergenceless vector field in two variables
can be presented in the form
$X=D_{1,2}(u)=\der_1(u)\der_2-\der_2(u)\der_1.$ We find that
$$s_5(D_{12}(u_1),D_{12}(u_2),D_{12}(u_3),D_{12}(u_4),D_{12}(u_5))=
-3 D_{12}([u_1,u_2,u_3,u_4,u_5]),
$$
where $[u_1,u_2,u_3,u_4,u_5]$ is the following determinant
$$
[u_1,u_2,u_3,u_4,u_5]=
\left\vert\begin{array}{ccccc}
\der_1u_1&\der_1u_2&\der_1u_3&\der_1u_4&\der_1u_5\\
\der_2u_1&\der_2u_2&\der_2u_3&\der_2u_4&\der_2u_5\\
\der_1^2u_1&\der_1^2u_2&\der_1^2u_3&\der_1^2u_4&\der_1^2u_5\\
\der_1\der_2u_1&\der_1\der_2u_2&\der_1\der_2u_3&\der_1\der_2u_4&
\der_1\der_2u_5\\
\der_2^2u_1&\der_2^2u_2&\der_2^2u_3&\der_2^2u_4&\der_2^2u_5\\
\end{array}\right\vert$$
For example,
$$s_5(\der_1,\der_2,x_1\der_2,x_1\der_1-x_2\der_2,x_2^2\der_1)=$$
$$(1/6)s_5(D_{12}(x_1),D_{12}(x_2), D_{12}(x_1^2),D_{12}(x_1x_2), D_{12}(x_2^3))=$$
$$-\frac{1}{2}D_{12}(
\left\vert\begin{array}{ccccc}
1&0&2x_1&x_2&0\\
0&1&0&x_1&3x_2^2\\
0&0&2&0&0\\
0&0&0&1&0\\
0&0&0&0&6x_2\\
\end{array}\right\vert)= -6\, D_{12}(x_2)=6 \der_1.$$

We can construct similar formulas for $n=2,3$ and $4.$
Formulas for $n=3,4$ are too
huge to present in an article.
In this  paper we present such a formula for $s_6.$
To calculate a $6$-commutator on $Vect(2)$ one needs to
calculate a sum of fourteen $6\times 6$ determinants (see
theorem~\ref{s6}.)
All properties given above for $5$-commutators can be generalized for
the $6$-commutator $s_6$ on $Vect(2).$

One more property that a $5$-commutator has, but a $6$-commutator does not,
is as follows. It is not true that a composition of adjoint derivations is a
derivation of Lie algebras, but
$$ad\,[X,Y]=[ad\,X,ad\,Y], \quad \forall X,Y\in Vect(n).$$
In general it is also false that $s_k(ad\,X_1,\ldots,ad\,X_k)$ is
a derivation. Nevertheless, for $Vect_0(2),$ $k=5,$ it is.
Moreover, $s_5(ad\,X_1,\ldots,ad\,X_5)$ is an adjoint
derivation and its potential can be constructed by the $5$-~commutator:
$$ad\,s_5(X_1,\ldots,X_5)=s_5(ad\,X_1,\ldots,ad\,X_5),$$
for any $X_1,\ldots,X_5\in Vect_0(2).$
A similar result for a $6$-commutator is no longer true.
For example,
$$F=s_6(ad\,\der_1,ad\,\der_2,ad\,x_1\der_1,ad\,x_2\der_1,ad\,x_1\der_2,
ad\,x_2\der_2)\in End\,Vect(2),$$
as a linear operator on $Vect(2)$ is defined by
$$(u_1\der_1+u_2\der_2)F=6(\der_1\der_2(u_1)+\der_2^2(u_2))\der_1-6(\der_2^2(u_1)+\der_1\der_2(u_2))\der_2.$$
We see that $F$ has nontrivial quadratic differential part, so it is not
even a derivation of the Lie algebra $(Vect(2),[\;,\;]).$

Note also the following relation between $5$ and $6$-commutators
and divergences of vector fields:
$$s_6^{rsym.r}(X_1,\ldots,X_6)=\sum_{i=1}^6(-1)^{i+1}Div\,X_i\, s_5^{rsym.r}(X_1,\ldots,\hat{X_i},\ldots,X_6),$$
for any $X_1,\ldots,X_6\in Vect(2).$ Here one can change on the right hand $s_5^{rsym.r}$ to $s_5,$ despite of the fact that $s_5$ is not well defined on $Vect(2).$
Quadratic differential part of a $5$-commutator can be presented as a sum of
three determinants (see lemma~\ref{quadratic part for s5}).
All quadratic differential terms of $s_5$ are cancelled in
taking alternative sum:
$$s_6^{rsym.r}(X_1,\ldots,X_6)=\sum_{i=1}^6(-1)^{i+1}Div\,X_i\,
s_5(X_1,\ldots,\hat{X_i},\ldots,X_6),$$
for any $X_1,\ldots,X_6\in Vect(2).$ Recall that $s_6=s_6^{rsym.r}$ on
$Vect(2).$ Notice that here $6$ and $5$ can not be changed to lower numbers.
For example,
$$s_5(\der_1,\der_2,x_1\der_1,x_2\der_1,x_1\der_2)-Div(x_1\der_1)s_4(\der_1,\der_2,x_2\der_1,x_2\der_2)=-3 \der_1\der_2\ne 0.$$

Notice that these results valid for $5$ and $6$-commutators do not valid
for lower degree commutators.
Namely $s_3,s_4,$ for $Vect_0(2)$ and
$s_3,s_4,s_5$ for $Vect(2)$ have no such properties.
One can state some weaker versions of these statements.
Let $gl_2+{\bf C}^2$ is a semi-direct sum of $gl_2$ and
standard $2$-dimensional module ${\bf C}^2.$
For example, if $k=3,4,5,$ then
$$[X,s_k^{rsym.r}(X_1,\ldots,X_k)]=\sum_{i=1}^k
s_k^{rsym.r}(X_1,\ldots,X_{i-1},[X,X_i],X_{i+1},\ldots,X_k),$$
for any $X\in  gl_2+{\bf C}^2,$ and
$X_1,\ldots,X_k\in Vect(2).$

Since composition of linear operators is a tensor operation,
$s_{N}$ for $N=n^2+2n-2$ is also a tensor operation. It is
a natural differential operation on the space of vector fields in
sense of \cite{Kirillov}. Identities of Witt algebra as a
Lie algebra were considered  by Razmyslov (see
\cite{Razmyslov}). One of his conjectures about minimality of Lie
identity of Lie algebra $W_n$ of degree $(n+1)^2$ is still open in
general case. For $n=2$ it was proved in \cite{Kagarmanov}.
Identities of Witt algebra as a right-symmetric algebra and their
connection with $N$-commutators was considered in \cite{Dzh00}.

\section{Preliminary facts and main definitions}

\subsection{${\mathcal D}$-differential algebras}
For simplicity, all vector spaces are considered over the field of
complex numbers ${\bf C}.$ By
 ${\bf Z}_+$ we denote the set of nonnegative integers.

A vector space $A$ is called a {\it $k$-algebra}
with multiplication $\omega$ and denoted
$A=(A,\omega)$ , if $\omega$ is a polylinear map
$A\times \cdots \times A\rightarrow A$
with $k$ arguments. Usually multiplications are understood as a bilinear
maps and instead of $\omega(a,b)$ one writes something like
$a\circ b$ or $a \cdot b.$ In such cases we will omit $2$ in the name
$2$-algebra and write $A=(A,\circ)$ or $A=(A,\cdot).$

Let $U$ be an associative commutative algebra with commuting system of
derivations ${\mathcal D}=\{\der_1,\ldots,\der_n\}.$
Call $U$ {\it $\mathcal D$-differential} algebra.
Let $L=U{\mathcal D}$ be its Lie algebra of derivations
$u\der_i: i=1,\ldots,n$ by the multiplication
$$[u\der_i,v\der_j]=u\der_i(v)\der_j-v\der_j(u)\der_i.$$
One can identify derivations as a vector fields and
consider a space of vector fields $Vect(n)=\langle u\der_i : u\in U, i=1,\ldots,n\rangle.$
So, $U{\mathcal D}=(Vect(n),[\;,\;]),$ i.e., as a linear spaces
$U{\mathcal D}$ and $Vect(n)$ are
coincide. Below we will endow $Vect(n)$ with other multiplications.
For example, $U{\mathcal D}^{rsym}=(Vect(n),\circ),$
the vector space $Vect(n)$ by the multiplication given by
$$u\der_i\circ v\der_j=v\der_j(u)\der_i,$$
is {\it right-symmetric}. i.e., it satisfies the right-symmetric identity
$$(X_1,X_2,X_3)=(X_1,X_3,X_2),$$
where $(X_1,X_2,X_3)=X_1\circ (X_2\circ X_3)-(X_1\circ X_2)\circ
X_3$ is an associator. Here we would like to pay attention to some
differences of our notations with usual ones. Usually action of
vector field to a function is denoted as $X(u),$ a vector field
$X$ is in the left hand and an argument $u$ is on the right.  But
in considering of right-symmetric algebras and its Lie algebras we
denote such actions as $(u)X.$ Therefore, the commutator given
above for $U{\mathcal D}$ as a Lie algebra and the commutator of a
Lie algebra obtained from right-symmetric algebra $U{\mathcal
D}^{rsym}$ are differ by signs.

\subsection{Witt algebra}
The basic example of ${\mathcal D}$-differential algebra
for us is $U={\bf C}[x_1,\ldots,x_n]$ with a system of
partial derivations $\der_i=\der_i/\der x_i, i=1,\ldots,n.$
In this case the algebra $U{\mathcal D}$ is denoted by $W_n.$
Call Lie algebra $W_n$ a {\it Witt algebra.}
Unless stated otherwise, by $U$ and $L=U{\mathcal D}$
we will understand these algebras.

Note that $U$ has a natural grading
$$U=\oplus_{s\ge 0}U_s,\quad  U_s U_k\subseteq U_{s+k},$$
$$U_s=\langle x^{\alpha}: \alpha\in \Gamma_n^+, |\alpha|=\sum_{i=1}^n\alpha_i=s\rangle, $$
and a filtration
$$U={\mathcal U}_{0}\supseteq {\mathcal U}_{1}\supseteq
{\mathcal U_{2}}\supseteq \cdots, \quad
{\mathcal U}_{s} {\mathcal U}_{k}\subseteq {\mathcal U}_{s+k},$$
$${\mathcal U}_{s}=\oplus_{s\ge k}U_k.$$

These grading and filtration induce a grading and a filtration on the
vector space of $L=W_n.$ Then a grading on $L$ is given by
$$L=\oplus_{s\ge -1} L_s, \quad
L_s=\langle x^{\alpha}\der_i: \alpha\in \Gamma_n^+, |\alpha|=s+1, i=1,\ldots,n\rangle $$
and a filtration
$$L={\mathcal L}_{-1}\supseteq {\mathcal L}_0\supseteq {\mathcal L}_1\supseteq
\cdots ,\quad {\mathcal L}_k=\oplus_{s\ge k} L_s.$$
These grading and filtration are compatible with the right-symmetric
multiplication:
$$L_s\circ L_k\subseteq L_{s+k},\quad
{\mathcal L}_s\circ {\mathcal L}_k\subseteq {\mathcal L}_{k+s},$$
for any $k,s\ge -1.$ In particular, they induce grading and
filtration of Witt algebra $W_n$ as a right-symmetric algebra and as a
Lie algebra.

For example, as vector spaces,
$$U_0=\langle1\rangle\cong {\bf C},$$
$$U_1=\langle x_i: i=1,\ldots n\rangle\cong {\bf C}^n,$$
$$U_2=\langle x_ix_j: i\le j, i,j=1,\ldots n\rangle\cong {\bf C}^{(n+1)n/2},$$
and
$$L_{-1}=\langle\der_1,\ldots,\der_n\rangle\cong {\bf C}^n,$$
$$L_0=\langle x_i\der_j: i,j=1,\ldots,n>  \cong {\bf C}^{n^2},$$
$$L_1=\langle x_ix_j\der_s: i,j,s=1,\ldots,n\rangle\cong {\bf C}^{(n+1)n^2/2}.$$

By the grading property of Witt algebra $W_n^{rsym}$ we see that
$$L_0\circ L_s \subseteq L_s$$
for any $s\ge -1.$ In particular, $L_0$ is a right-symmetric subalgebra
of $L=W_n^{rsym}.$ Hence it is a Lie subalgebra of $W_n.$
One can see that the subalgebra $L_0\subset W_n^{rsym}$
is not only right-symmetric. It is associative and isomorphic to
the matrix algebra $Mat_n.$ It is easy to see that $L_0$ as a Lie algebra
is isomorphic to $gl_n;$ moreover,
$$[L_0,L_0]=\langle x_j\der_s, x_i\der_i-x_{i+1}\der_{i+1}: 1\le i<n, 1\le j,s\le n,
j,\ne s,\rangle\cong sl_2.$$
$$L_0=[L_0,L_0]\oplus \langle e\rangle, \quad e=\sum_{i=1}^nx_i\der_i.$$
So, we have $sl_n$-module structures on $L_k.$ Let $\pi_1,\ldots,
\pi_{n-1}$ be fundamental weights of $sl_n$ and $R(\gamma)$ denote
the irreducible module over Lie algebra $sl_n$ with the highest
weight $\gamma.$

For example, as a $sl_n$-modules,
$$L_{-1}\cong R(\pi_{n-1}),$$
$$L_0\cong R(\pi_1+\pi_{n-1})\oplus R(0),$$
$$L_1\cong R(2\pi_1+\pi_{n-1})\oplus R(\pi_1).$$

Let
$$Div: L\rightarrow U, \quad
X=\sum_{i=1}^nu_i\der_i\mapsto \sum_{i=1}^n\der_i(u_i),
$$
be a divergence map. Then
$$Div\,[X,Y]= X (Div(Y))-Y (Div(X)),$$
for any $X,Y,Z\in L.$ In particular,
$$Div\,[a,X]=a(Div(X)),$$
for any $a\in L_0$ and $X\in L.$ So, we have an
exact sequence of $sl_n$-modules
$$0\rightarrow {\bar L}_k\rightarrow L_k\stackrel{Div}
{\rightarrow} U_{k-1}\rightarrow 0,$$
where
$${\bar L}_k=\langle X\in L_k: Div\,X=0\rangle.$$
In other words, the following isomorphisms of $sl_n$-modules holds
$$L_k\cong {\bar L}_k\oplus{\tilde L}_k,$$
$${\bar L}_k\cong R((k+1)\pi_1+\pi_{n-1}),\quad \tilde L_k\cong R(k\pi_1).$$

\subsection{Root system in Witt algebra}
Let $\mathcal E$ be a subset of ${\bf Z}^n$
consisting of elements of the form $\alpha-\epsilon_i,$ where
$\alpha\in {\bf Z}^n_+$ and $\epsilon_i=(0,\ldots,0,1,0,\ldots,0)$ ($1$ is in
the $i$-th place), $i=1,\ldots,n.$

Note that $L=W_n$ has a Cartan subalgebra $H=\langle h_1,\ldots,h_n\rangle $ where
$h_i=x_i\der_i, i=1,\ldots,n$ and $W_n$ has a root decomposition
with the root system $\mathcal E:$
$$L=\oplus_{\alpha\in \mathcal E}L_{\alpha},\quad
[L_{\alpha},L_{\beta}]\subseteq L_{\alpha+\beta},$$
$$L_{\alpha}=\langle x^{\alpha+\epsilon_j}\der_j: \alpha+\epsilon_j\in {\bf Z}^n_+,
j=1,\ldots,n\rangle, \quad \alpha\in \mathcal E.$$ Moreover, the
actions of $h_i$ generate  a grading on $L$-module $U:$
$$U=\oplus_{\alpha\in {\bf Z}^n_+}U_{\alpha},$$
$$U_{\alpha} U_{\beta}\subseteq U_{\alpha+\beta}, \forall \alpha,\beta\in {\bf Z}^n_+,$$
$$U_{\alpha} L_{\beta} \subseteq U_{\alpha+\beta}, \quad \forall \alpha\in
{\bf Z}^n_+, \beta\in \mathcal E.$$
So, we have constructed an $\mathcal E$-grading on $W_n$ and on its
$(L,U)$-module $U.$ Here we suppose that ${\bf Z}^n_+\subseteq \mathcal E$
and $L_{\alpha}=0$ if $\alpha\not\in \mathcal E$ and $U_{\alpha}=0$ if
$\alpha\not\in {\bf Z}^n_+.$

Let $L=W_n$ and assume that all modules $M_1\ldots, M_{k+1}$ coincide with the
adjoint module $L.$ Let $\psi\in T^k(L,L)^{L_{-1}}:$
$$\der_i(\psi(X_1,\ldots,X_k))=\sum_{i=s}^k\psi(X_1,\ldots,X_{s-1},
\der_i(X_s),X_{s+1},\ldots,X_k),$$ for any $X_1,\ldots,X_k\in L.$
The notation $T^k$  means here that we consider polylinear maps with
$k$ arguments, but arguments are not necessary skew-symmetric.
Suppose that $\psi$ is $\mathcal E$-graded:
$$\psi(L_{\alpha^1},\ldots,L_{\alpha^k})\subseteq
L_{\alpha^1+\cdots+\alpha^k},$$
for any $\alpha^1,\ldots,\alpha^k\in \mathcal E.$

\subsection{${\mathcal D}$-invariant  polynomials }

Let $f=f(t_1,\ldots,t_k)$ be an associative noncommutative polynomial,
$f=\sum_{i_1,\ldots,i_l}\lambda_{i_1,\ldots,i_l}t_{i_1}\ldots t_{i_l}.$
I.e., $f$ is a linear combination of monomials $t_{(i)}=t_{i_1}\ldots t_{i_l},$
where $i_1,\ldots,i_l$ run through elements $1,\ldots,k,$ possibly repeated.
Later on we will make substitutions of variables $t_i$ by elements
of some $2$-algebras. Since our algebras may be non-associative
in such cases (unless stated otherwise)
we will suppose that any monomial $t_{(i)}$ has a {\it left-normed} bracketing,
i.e., $t_{(i)}=(\cdots((t_{i_1}t_{i_2})t_{i_3})\cdots)t_{i_l}.$

Let $A=(A,\circ)$ be an algebra with multiplication $\circ.$
Suppose that $B$ is a  subspace of $A.$
Since $f$ is a linear combination of monomials of the form
$t_{(i)},$ one can substitute instead of $t_i$
elements of $B$ and calculate $T$ using the multiplication of
the algebra $A.$ So, we obtain a map
$f_B: B\times \cdots\times B\rightarrow A$ defined by
$$f_B(b_1,\ldots,b_k)=f(b_1,\ldots,b_k).$$
Sometimes we endow $A$ by several multiplications. In such cases,
we will write $f_{B}^{\circ}$ instead of $f_B,$ if $A$ is considered
as an algebra with multiplication $\circ.$ Notice that $B$ may be
not closed under multiplication $\circ.$
Whenever it is clear from the context to what $B$ we make substitutions,
we reduce the notation $f_B$ to $f.$

Suppose now that $U$ is a ${\mathcal D}$-differential algebra. We
endow the space of differential operators $Dif\!f(n)$ by three
multiplications. Namely, the multiplications $\cdot,$ $\circ$ and
$[\;,\;]$  stand for composition, right-symmetric multiplication
and commutators. We will  sometimes write $f_{Vect(n)}^{\cdot}=f,
f_{Vect(n)}^{\circ}=f^{rsym.r}.$

Call a polynomial $f=f(t_1,\ldots,t_k)$
or more precisely $f_{Dif\!f(n)},$ is ${\mathcal D}$-invariant, if
$$[f(X_1,\ldots,X_k),\der_i]=\sum_{s=1}^kf(X_1,\ldots,X_{s-1},[X_s,\der_i],X_{s+1},\ldots,X_k)$$
for any $X_1,\ldots,X_k\in Dif\!f(n)$ and $i=1,\ldots,n.$ Here we do not
specify what kind of multiplications in $Dif\!f(n)$ we use
in the calculation of $f_{Dif\!f(n)}$. The reason is
that $$[X,\der_i]=X\cdot \der_i-\der_i\cdot X=X\circ \der_i-\der_i\circ X,$$
for any $X\in Dif\!f(n).$ By this observation, $\mathcal D$-invariance of
$f_{Dif\!f(n)}^{\cdot}$ and $f_{Dif\!f(n)}^{circ}$ are equivalent notions.

\subsection{Primitive commutators}

Assume that $g=g(t_1,\ldots,t_k)$ is a skew-symmetric
multilinear polynomial with $k$ arguments. We call $g$ a
{\it $k$-commutator} on some class of vector
fields,  if for any $k$ vector fields $X_1,\ldots,X_k$
of this class $g(X_1,\ldots,X_k)$ is again a vector field
of this class.

Suppose that $f$ is some polynomial with left-normed brackets.
Let $(A,[\;,\;])$ be a Lie algebra with vector space $A$ and
commutator $[\;,\;].$ As we have explained above,
$f_A^{[\;,\;]}: A\times \cdots
\times A\rightarrow A$ is a map with $k$ arguments that is obtained
from $f$ by substituting elements of algebra $A$
instead of parameters $t_j$ and using the commutator $[\;,\;].$

Suppose now that $(A, [\;,\;])$ is a Lie algebra for some class of
vector fields.  Then $f_A^{[\;,\;]}$ is the standard  $k$-commutator
for any vector field algebra $A.$ Call $k$-commutator $g$
{\it primitive} on $Vect(n)$
if $g_{(Vect(n)}^{\cdot}$ can not be represented as
$f_{Vect(n)}^{[\;,\;]},$ for any
left-normed polynomial $f.$
It is easy to give a similar definition for $Vect_0(n).$

The right-symmetric multiplication on $Vect(n)$  is given by
$$u\der_i\circ v\der_j=v\der_j(u)\der_i.$$
This multiplication satisfies the {\it right-symmetric identity}
$$(X,Y,Z)=(X,Z,Y),$$
where $(X,Y,Z)=X\circ (Y\circ Z)-(X\circ Y)\circ Z$ is the
associator of $X,Y,Z.$ In other words,
$$X\circ [Y,Z]=(X\circ Y)\circ Z-(X\circ Z)\circ Y,$$
for any $X,Y,Z\in Vect(n).$ So, we obtain a right-symmetric
algebra $(Vect(n),\circ).$
Any right-symmetric algebra is Lie-admissible. In particular,
$(Vect(n), [\;,\;]),$ i.e., the space $Vect(n)$ with the commutator
$[X,Y]=X\circ Y-Y\circ X$ is a Lie algebra. As we have noticed, the
commutator $X\circ Y-Y\circ X$ coincides here with the standard
commutator $X \cdot Y-Y\cdot X.$

\subsection{$n$-Lie algebras, $(n-1)$-left commutative algebras and $sh$-algebras}
For vector spaces $M$ and $N$ denote by $T^k(M,N)=Hom(M^{\otimes k},N)$ the
space of polylinear maps with $k$ arguments $\psi: M\times \cdots \times M\rightarrow N$ if $k>0,$ and set $T^0(M,N)=N$ and $T^k(M,N)=0$ if $k<0.$
Let $T^*(M,N)=
\oplus_kT^k(M,N).$ Let $\wedge^kM$ be the $k$-th exterior power of $M$ and $C^k(M,N)=Hom(\wedge^kM,N)\subseteq T^k(M,N)$ be a space of skew-symmetric
polylinear maps with $k$ arguments in $M$ and $C^0(M,N)=N,$ $C^k(M,N)=0,$ if $k<0.$ Set
$C^*(M,N)=\oplus_kC^k(M,N).$

Let $\Omega=\{\omega_1,\omega_2,\ldots\}$ be a set of polylinear maps
$\omega_i\in C^i(A,A).$ Let  $(A,\Omega)$ be an algebra with vector space $A$ and signature $\Omega$ \cite{Kurosh}.
This means that $A$ is endowed by the $i$-ary multiplication $\omega_i.$ Call $A$ an {\it $\Omega$-algebra} and in case where $\Omega=\{\omega_i\}$
set $A=(A,\omega_i).$ For $\omega\in \Omega,$ write $|\omega|=i$
if $\omega\in C^i(A,A).$

Suppose $(A,\omega)$ is an algebra.
It is called an {\it $n$-Lie} \cite{Filipov1} or a {\it Nambu} algebra,
if
$$
\omega(a_1,\ldots,a_{n-1},\omega(a_{n},\ldots,a_{2n-1}))=$$
$$\sum_{i=0}^{n-1}\omega(a_n,\ldots,\omega(a_1,\ldots,a_{n-1},a_{n+i}),a_{i+1},\ldots,a_{2n-1}),$$
{\it $(n-1)$-left commutative,} if
$$\sum_{\sigma\in Sym_{2n-2}}sign\,\sigma\, \omega(a_{\sigma(1)},\ldots\omega_{\sigma(n-1)},\omega(a_{\sigma(n)},\ldots,a_{\sigma(2n-2)},a_{2n-1}))=0,$$
and {\it $n$-homotopical Lie}, if
$$\sum_{\sigma\in Sym_{2n-1}}sign\,\sigma\, \omega(a_{\sigma(1)},\ldots\omega_{\sigma(n-1)},\omega(a_{\sigma(n)},\ldots,a_{\sigma(2n-2)},a_{\sigma(2n-1)}))=0,$$ for any $a_1,\ldots,a_{2n-1}\in A.$

Finally, an algebra $(A,\Omega),$ where $\Omega=\{\omega_1,\omega_2,\ldots\}$
is called {\it homotopical Lie} or $sh$-Lie \cite{StashefLada}, if
$$\sum_{i+j=k+1, i,j\ge 1}(-1)^{(j-1)i}sign\,\sigma\,
\omega_j(\omega_i(a_{\sigma(1)},\ldots,a_{\sigma(i)}),
a_{\sigma(i+1)},\ldots,a_{\sigma(i+j-1)})=0,$$
for any $k=1,2,\ldots,$ and for any $a_1,\ldots,a_{i+j-1}\in A.$ In
particular, an $n$-homotopical Lie algebra is an $sh$-algebra if $\Omega$
is concentrated in
only one non-zero multiplication $\omega_n,$ i.e., $\omega_i=0,$ if $i\ne n$.
Suppose now that $\Omega$ is concentrated in two elements
$\omega_2$ and $\omega_n.$ Then  the condition that $(A,\Omega)$ is $sh$-Lie
means that $(A,\omega_2)$ is a Lie algebra,
$(A,\omega_n)$ is a $n$-homotopical Lie and $\omega_n$ is a $n$-cocycle of
the adjoint module of the Lie algebra $(A,\omega_2).$

In \cite{Dzhumavronskian} it is established that, over the field of
characteristic 0, any $n$-Lie algebra is $(n-1)$-left commutative and
any $(n-1)$-left commutative algebra is $n$-homotopical Lie. We
prove that $(Vect_0(2),s_5)$ and $(Vect(2),s_6)$ are
correspondingly {\it $4$- and $5$-left commutative}. Hence {\it
$(Vect_0(2),s_5)$ is $5$-homotopical Lie} and {\it $(Vect(2),s_6)$
is $6$-homotopical Lie.} But they are not, correspondingly, 5- and
6-Lie.

\section{Main results}
\begin{thm}\label{N-commutator} Let $N=n^2+2n-2.$ Then

{\rm (i)} $s_k=s_k^{rsym.r}$ on $Vect(n),$ for any $k\ge N.$
In particular, $s_k$ is well defined on $Vect(n),$ if $k\ge N.$

{\rm (ii)} $s_k=0$ is an identity on $Vect(n)$ for any  $k> N+1.$

{\rm (iii)} $(Vect(n),s_N)$ is $(N-1)$-left commutative.
In particular, $(Vect(n),s_N)$ is $N$-homotopical Lie.

{\rm (iv)} $ad\,X\in Der\,(Vect(n), s_k)$ for any $X\in Vect(n)$ and for any
$k\ge N.$

{\rm (v)} $ad\,X\in Der\,(Vect(n),s_k)$ for any $X\in Vect(n)$ such that
$\der_i\der_j(X)=0, i,j=1,\ldots,n.$ Here $k$ is any integer $>0.$
\end{thm}

\begin{thm} \label{5-commutator} {\rm (i)} $s_5\ne 0$ on $Vect_0(2).$

{\rm (ii)} $s_5$ is a $5$-commutator on $Vect_0(2).$

{\rm (iii)} $s_6=0$ is an identity on $Vect_0(2).$

{\rm (iv)} $5$-commutator $s_5$ on $Vect_0(2)$ is primitive.

{\rm (v)} $s_5=s_5^{rsym.r}$ on $Vect_0(2).$

{\rm (vi)} $ad\,s_5(X_1,\ldots,X_5)= s_5(ad\,X_1,\ldots,ad\,X_5)$ for
any $X_1,\ldots,X_5\in Vect_0(2).$

{\rm (vii)} $(Vect_0(2),s_5)$ is a $4$-left commutative algebra.

{\rm (viii)} $(Vect_0(2),\{s_2,s_5\})$ is an $sh$-Lie algebra.
\end{thm}

\begin{thm}\label{6-commutator}
{\rm (i)} $s_6\ne 0$ on $Vect(2).$

{\rm (ii)} $s_6$ is a $6$-commutator on $Vect(2).$

{\rm (iii)} $s_7=0$ is an identity on $Vect(2).$

{\rm (iv)} $6$-commutator $s_6$ on $Vect(2)$ is primitive.

{\rm (v)} $s_6=s_6^{rsym.r}$ on $Vect(2).$

{\rm (vi)} For any $X_1,\ldots,X_6\in Vect(2),$
$$s_6(X_1,\ldots,X_6)=\sum_{i=1}^6(-1)^{i+1}Div\,X_i\,s_5(X_1,\ldots,\hat{X_i},\ldots,X_6).$$

{\rm (vii)} $(Vect(2),s_6)$ is a $5$-left commutative algebra.

{\rm (viii)} $(Vect(2),\{s_2,s_6\})$ is an $sh$-Lie algebra.
\end{thm}

A natural question arises: Is it possible to construct nontrivial
$N$-~commutators on $Vect(n)$ for $n>2$ ?
If $N=n^2+2n-2,$ then $s_N$ is a  well defined
$N$-commutator on $Vect(n), n>1$  (theorem~\ref{N-commutator}, (i)).
One can prove that $s_{N-1}$ is also a well defined operation on $Vect_0(n).$

{\bf Conjecture}. \begin{itemize}
\item $S_N(X_1,\ldots,X_N)\ne 0$ for some $X_1,\ldots,X_N\in
Vect(n).$
\item $S_{N-1}(X_1,\ldots,X_{N-1})\ne 0,$ for some $X_1,\ldots,X_N\in
Vect_0(n).$
\end{itemize}

We have checked this conjecture by a computer  for $n=2, 3 ,4. $

\section{$k$-commutators by right-symmetric multiplications}

The aim of this section is to prove
that for any $n$ there exists
$N=N(n)$ such that for any $k\ge N,$ $s_k=s_k^{rsym.r}$ on $Vect(n).$

For $X=\sum_{\alpha\in{\bf Z}^n-+}  u_{\alpha}\der^{\alpha}\in Dif\!f(n),$ where $u_{\alpha}\in U,$ set $\vert X\vert =s$
if $u_{\alpha}\ne 0$ for some
$\alpha\in {\bf Z}^n_+,$ $\vert\alpha\vert=s$ but $u_{\beta}=0$ for any
$\beta\in {\bf Z}^n_+$ such that $\vert \beta \vert<s.$
Set also $||X||=q$ if $u_{\alpha}\ne 0$ for some $\alpha\in {\bf Z}^n_+,$
$|\alpha|=q,$ but $u_{\beta}=0$ for all $\beta\in {\bf Z}^n_+$ such that
$|\beta|>q.$ Notice that
$$|X+Y|\ge min\{|X|,|Y|\},$$
$$||X+Y||\le  max\{||X||,||Y||\},$$
for any $X,Y\in Dif\!f(n).$

Define a multiplication  $Dif\!f(n)\times Dif\!f(n)\rightarrow Dif\!f(n),
(X,Y)\mapsto X\bullet Y$ by
$$u\der^{\alpha}\bullet v\der^{\beta}=
\sum_{\gamma\in {\bf Z}^n_+, \gamma\ne 0}{\beta\choose\gamma}v\der^{\beta-\gamma}(u)\der^{\alpha+\gamma},$$
for any $\alpha,\beta\in {\bf Z}^n_+.$
Prolong right-symmetric multiplication $\circ$ from a space of differential
operators of first order to a space of differential operators of any order by
$$u\der^{\alpha} \circ v\der^{\beta}=v\der^{\beta}(u)\der^{\alpha}.$$
We see that
$$X\cdot Y=X\circ Y+X\bullet Y,$$
$$|X\bullet Y|>|X|, \mbox{\;if\;} X\bullet Y\ne 0,$$
for any $X,Y\in Dif\!f(n).$

\begin{lm} Suppose that
$X_1,\ldots,X_k\in Dif\!f(n),$ are linear differential operators such that
$\vert X_i\vert\ge s$ for any $i=1,2,\ldots,k$ and
$\vert\vert s_k(X_1,\ldots,X_k)\vert\vert\le s.$
Then
$$s_k(X_1,\ldots,X_k)=pr_s(s_k^{rsym.r}(X_1,\ldots,X_k)),$$
where $pr_s: Dif\!f(n)\rightarrow Dif\!f(n)$ a projection map to
$\langle u\der^{\alpha}: \alpha\in {\bf Z}^n_+, |\alpha|=s\rangle.$
\end{lm}

{\bf Proof.}
By the formula $X\cdot Y=X\circ Y+X\bullet Y$ one can write a composition
$X_{\sigma(1)} \cdot \ldots \cdot X_{\sigma(k)}$ as a
linear combination of elements of the form
$X_{\sigma}'=
(\ldots ((X_{\sigma(1)} \circ X_{\sigma(2)}) \ldots )\circ X_{\sigma(k)},$
where no $\bullet$ is appeared, and elements of the form
$X''_{\sigma}=
(\ldots ((X_{\sigma(1)} \star X_{\sigma(2)}) \ldots )\star X_{\sigma(k)},$
where $\star= \circ$ or $\bullet$ and the number of the
$\bullet$'s is at least one.
Notice that $|X\bullet Y|> |X|$ if $|X\bullet Y|\ne 0.$ Therefore,
$$|X_{\sigma}''|>s,$$
if $X_{\sigma}''\ne 0.$
Since
$$s_k(X_1,\ldots,X_k)=\sum_{\sigma\in Sym_k}X_{\sigma}'+X_{\sigma}'',$$
$$||s_k(X_1,\ldots,X_k)||\le s,$$
$$|\sum_{\sigma\in Sym_k}sign\,\sigma\,X_{\sigma}'|\ge s,$$
$$|\sum_{\sigma\in Sym_k}sign\,\sigma\,X_{\sigma}''|>s,$$
which means that
$$||s_k(X_1,\ldots,X_k)||=|s_k(X_1,\ldots,X_k)|=s$$
and
$$s_k(X_1,\ldots,X_k)=pr_s(s_k^{rsym.r}(X_1,\ldots,X_k)).$$

\begin{crl} \label{14fevr}
Suppose that $s_k(X_1,\ldots,X_k)\in Vect(n)$ for
$X_1,\ldots,X_k\in Vect(n).$
Then
$$s_k(X_1,\ldots,X_k)=s_k^{rsym.r}(X_1,\ldots,X_k).$$
\end{crl}

\begin{lm} Suppose that
$s_k=s_k^{rsym.r}$ on $Vect(n)$ and
$D$ is a derivation of $(Dif\!f(n),\cdot)$ that saves subspace $Vect(n).$
Then $D$ is a derivation of $k$-algebra $(Vect(n),s_k^{rsym.r}),$ i.e.,
$$D(s_k^{rsym.r}(X_1,\ldots,X_k))=\sum_{i=1}^ks_k^{rsym.r}(X_1,\ldots,X_{i-1},D(X_i),X_{i+1},\ldots,X_k),$$
for any $X_1,\ldots,X_k\in Vect(n).$
\end{lm}

{\bf Proof.} We have
$$D(s_k^{rsym.r}(X_1,\ldots,X_k))=$$
(corollary \ref{14fevr})
$$D(s_k(X_1,\ldots,X_k))=$$
(since $D\in Der\,Dif\!f(n)$)
$$=\sum_{i=1}^ks_k(X_1,\ldots,X_{i-1},D(X_i),X_{i+1},\ldots,X_k).$$
Since $D(Vect(n))\subseteq Vect(n)$ by our condition,
$$s_k(X_1,\ldots,X_{i-1},D(X_i),X_{i+1},\ldots,X_k)\in Vect(n),$$
for any $X_1,\ldots,X_k\in Vect(n).$  Thus by corollary \ref{14fevr},
$$s_k(X_1,\ldots,X_{i-1},D(X_i),X_{i+1},\ldots,X_k)=$$
$$s_k^{rsym.r}(X_1,\ldots,X_{i-1},D(X_i),X_{i+1},\ldots,X_k).$$
Hence,
$$D(s_k^{rsym.r}(X_1,\ldots,X_k))=
\sum_{i=1}^ks_k^{rsym.r}(X_1,\ldots,X_{i-1},D(X_i),X_{i+1},\ldots,X_k),$$
for any $X_1,\ldots,X_k\in Vect(n).$

\begin{crl} \label{alia}
Suppose that $s_k=s_k^{rsym.r}$ on $Vect(n).$ Then for any $X\in
Vect(n)$ the adjoint derivation $ad\,X$ generates also a derivation
of $(Vect(n),s_k^{rsym.r}).$ In particular, $ad\,X$ is
a derivation of the algebra $(Vect(n),s_{n^2+2n-2})$ for any $X\in
Vect(n).$  Similarly, $ad\,X$ is a derivation of
$(Vect_0(2),s_5)$ for any $X\in Vect_0(2).$
\end{crl}

\section{How to calculate $\mathcal D$-invariant forms}

\subsection{$(L,U)$-modules}
Let $U={\bf C}[x_1,\ldots,x_n]$ and $L=W_n=Der\,U.$

Let $M$ be $L$-module. Then the subspace of
$L_{-1}$-invariants, $M_0=M^{L_{-1}}=
\langle m\in M: (m)\der_i=0\rangle ,$ has a natural structure of
an $L_0$-module.
Make $M_0$ an $\mathcal L_0$-module by the trivial action of
${\mathcal L}_1.$
Call the ${\mathcal L}_1$-module $M_0$ the {\it base} of $L$-module $M.$

Let $M$ be a module over Lie algebra $L=W_n$. Call $M$ an
{it $(L,U)$-module,} if there exists an additional structure of $M$ as
a right $U$-module  such that
$$(mu)X=m[u,X]+(mX)u,$$
for any $m\in M, u\in U, X\in L.$

Let $M$ be an $(L,U)$-module and $M_0=M^{L_{-1}}$ is a base of $M$ as
an $L$-module.
Call $M$ an {\it $(L,U)$-module with base $M_0,$} if $M$ is an
$(L,U)$-module and
$M$ as a $U$-module is a free module with base $M_0.$

The main construction of $(L,U)$-modules (\cite{Dzhumavestnik})
is the following.

Let $L=W_n$ and $U$ be a natural $L$-module. Let $M_0$ be
an ${\mathcal L}_0$-module such that $M_0{\mathcal L}_1=0.$
Recall that $\mathcal L_0=\oplus_{s\ge 0}L_i=L_0\oplus {\mathcal L}_1.$
One can make ${\mathcal M}=U\otimes M_0$ an
$(L,U)$-module by
$$(u\otimes m)X=(u)X\otimes m+u\otimes \sum_{a\in V({\mathcal L}_0}E_a(X)u\otimes (m)a,$$
$$(u\otimes m)v=(v u)\otimes m,$$
where $V({\mathcal L})_0=\{x^{\alpha}\der_i: \alpha\in {\bf Z}^n_+,
i=1,\ldots,n\}$ is the standard basis of ${\mathcal L}_0.$ Recall that for
$L=W_n$ and $L=S_{n-1}=\langle X\in W_n: Div\,X=0\rangle $ the map $E_a: L\rightarrow U$
is defined by
$$E_{x^{\alpha}\der_j}
(v\der_j)=\delta_{i,j}\frac{\der^{\alpha}(v)}{\alpha !}.$$
Then ${\mathcal M}=U\otimes M_0$ will be an $(L,U)$-module with base $M_0.$
It has a standard grading induced  by a grading of $U:$
$${\mathcal M}=\oplus_{s\ge -q} {\mathcal M}_[s], \quad \M_{[s]}=\langle x^{\alpha}\otimes m: |\alpha|=\sum_{i=1}^n\alpha_i=s+q\rangle, $$
where one can assume that $M_0\cong \langle 1\otimes m: m\in M_0\rangle $ has the grading
degree $-q$ and $|x^{\alpha}\otimes m|=|\alpha|-q.$
For example,
$${\mathcal M}_{[2+q]}=\langle x_ix_j\otimes m: i\le j, m\in M_0\rangle,\quad \dim\,M_{[2+q]}={n+1\choose 2}
\dim\,M_0.$$

{\bf Example.} $U$  under the natural actions of $L=W_n$ and
$U$ is an $(L,U)$-module with base $\langle 1\rangle .$
The adjoint module of $L=W_n$ under the natural action
of $U,$ $(u\der_i,v)\mapsto u v \der_i,$ is also an $(L,U)$-module with base
$L_{-1}.$  This module is isomorphic to the $(L,U)$-module $U\otimes L_{-1}$
with the module structures given above.
If $L=S_{n-1}=\langle X\in W_n: Div\,X=0\rangle, $ then the adjoint module
has no structure of an $(L,U)$-module.

Let $ \mathcal M_1,\ldots,\M_k$ and $\mathcal M_{k+1}$ are
$(L,U)$-modules with bases $M_{1,0},\ldots,M_{k,0}$ and $M_{k+1,0},$
correspondingly. Let
$C^k({\mathcal M_1,\ldots,\M_k; \M_{k+1}})$ be a space of polylinear maps
$\psi: {\mathcal M_1\times \cdots\times \M_k\rightarrow \M_{k+1}}.$
Endow $C^k({\mathcal M_1,\ldots,\M_k;\M_{k+1}})$ with a structure of
$(L,U)$-module. The action of $L$ is natural:
$$(m_1,\ldots,m_k)(\psi X)=$$
$$((m_1,\ldots,m_k)\psi)X-
\sum_{i=1}^k\psi(m_1,\ldots,m_{i-1},(m_i)X,m_{i+1},\ldots,m_k).$$
The action of $U$ is given by
$$(m_1,\ldots,m_k)(\psi u)=((m_1,\ldots,m_k)\psi)u.$$
Here $\psi\in C^k({\mathcal M_1,\ldots,\M_k;\M})$ and $m_i\in \M_i, 1\le i\le k,$ $X\in L,$
$u\in U.$
It is easy to see that
$$(\psi u) X=\psi(u X)+(\psi X)u,$$
for any $\psi\in C^k({\mathcal M_1,\ldots,\M_k;\M_{k+1}}), X\in L$
and $u\in U.$ So, $C^k({\mathcal M_1,\ldots,\M_k;\M_{k+1}})$ is a
$(L,U)$-module. The algorithm how to find its base or how to
calculate {\it ${\mathcal D}$-invariant forms} is given in
\cite{Dzhumavestnik}.

We need results of \cite{Dzhumavestnik} in one particular case. Namely, we
use this algorithm in the calculation of $k$-commutators. This
algorithm is given below.

\subsection{Escorts and supports}

Let $\mathcal E$ be a root system on $W_n.$ associate with $\psi$ a
polylinear map $esc(\psi)\in T^k(L,L_{-1}),$ called the {\it
escort} of $\psi,$ by the rule
$$esc(\psi) (X_1,\ldots,X_k)=\psi(X_1,\ldots,X_k),$$
if $X_1\in L_{\alpha^1},\ldots, X_k\in L_{\alpha^k}$
and $\alpha^1+\cdots+\alpha^k$ has a form
$-\epsilon_s$ for some $s=1,\ldots,n.$ Here $\alpha^1,\ldots,\alpha^k$
are some roots from $\mathcal E.$ If $\alpha^1+\cdots+\alpha^k$ can not be presented in
the form $-\epsilon_s$ for some $s\in \{1,2,\ldots,n\}$ then, by definition,
$$esc(\psi)(X_1,\ldots,X_k)=0.$$
So, we have defined $esc(\psi)(X_1,\ldots,X_k)$ for root elements
$X_1,\ldots,X_k.$ Define now $esc(\psi)(X_1,\ldots,X_k)$ by
polylinearity for any $X_1,\ldots,X_k.$ Notice that $esc(\psi)$ is
correctly defined, i.e.,
$$\forall X_1,\ldots,X_k\in L \Rightarrow esc(\psi)(X_1,\ldots,X_k)\in
L_{-1}.$$

Notice also that $-\epsilon_s$ can not be presented as a sum of $k$ roots
in infinitely many ways. Therefore, the following subspace
$$supp_s(\psi)= \oplus_{\alpha^1,\ldots,\alpha^k\in \mathcal E,
\alpha^1+\ldots+\alpha^k=-\epsilon_s}^{} L_{\alpha^1}\otimes \cdots
\otimes L_{\alpha^k},$$
called  the {\it $s$-support of $\psi$},
is finite-dimensional for any $s\in \{1,\ldots,n\}.$
Call a finite-dimensional subspace
$$supp(\psi)=\oplus_{s=1}^n supp_s(\psi)$$
the {\it support of $\psi.$} So, escort of any ${\mathcal
D}$-invariant ${\mathcal E}$-graded map $\psi\in T^k(L,L)$ is
uniquely defined by the restriction to its support $supp(\psi).$

Let $V(L)$ be the standard basis of $L=W_n$ consisting of vectors
of the form $x^{\alpha}\der_i$ where $\alpha\in {\bf Z}^n_+$ and
$i\in\{1,\ldots,n\}.$
Denote by $V(\psi)$ a basis of $supp(\psi)$ which is collected by tensoring of
basic vectors of $V(L).$ We will write $(a_1,\ldots,a_k)$ instead of
$a_1\otimes \cdots\otimes a_k.$ So,
$$V(\psi)=\cup_{s=1}^nV_s(\psi),$$
where
$$V_s(\psi)=\{(a_1,\ldots,a_k):
a_l\in V(L_{\alpha^l}), \alpha^l\in{\mathcal E,}
\alpha^1+\ldots+\alpha^k=-\epsilon_s,
l=1,\ldots,k\}.$$

As was shown in \cite{Dzhumavestnik}, any $\mathcal E$-graded
${\mathcal D}$-invariant map can be uniquely restored by its escort.
Namely,
\begin{equation}
\label{escort}
\psi(X_1,\ldots,X_k)=\sum_{(a_1,\ldots,a_k)\in V(\psi)}
E_{a_1}(X_1)\cdots E_{a_k}(X_k)esc(\psi)(a_1,\ldots,a_k),
\end{equation}
where
$$E_{x^{\alpha}\der_i}(v\der_j)=\delta_{i,j}\frac{\der^{\alpha}(v)}{\alpha!}.$$

Denote by  $s_k$ the standard skew-symmetric associative polynomial
$$s_k=s_k(t_1,\ldots,t_k)=\sum_{\sigma\in Sym_k}sign\,\sigma\,t_{\sigma(1)}\cdots t_{\sigma(k)}.$$

\subsection{Cup-products}
Given an algebra $A$ with multiplication $\star$,
define the cup-product of cochain complexes $C^*(A,A)\times C^*(A,A)$
corresponding to a multiplication on $A.$ If $\psi\in C^k(A,A)$ and
$\phi\in C^l(A,A)$ then define $\psi\smile\phi\in C^{k+l}(A,A)$ by
$$\psi\smile \phi(a_1,\ldots,a_{k+l})=
\sum_{\sigma\in Sym_{k,l}}sign\,\sigma\,\psi(a_{\sigma(1)},\ldots,a_{\sigma(k)})\star \phi(a_{\sigma(k+1)},\ldots,a_{\sigma(k+l)}).$$
Here
$$Sym_{k,l}=\{\sigma\in Sym_{k+l}: \sigma(1)<\cdots \sigma(k), \sigma(k+1)<\cdots<\sigma(k+l)\}$$
is a set of $(k,l)$-shuffle permutations.

Suppose that $A$ has an associative  multiplication $\cdot$ and a
right-symmetric multiplication $\circ.$
Denote by $\overset{\cdot}\smile$ and $\overset{\circ}{\smile}$
cup-products induced by multiplications $\cdot$ and $\circ$ correspondingly.

If $B$ is a subspace of $A$ then one can consider cup-products
$$\smile: C^*(B,A)\times C^*(B,A)\rightarrow C^*(B,A).$$

We use the cup-products for $A=Dif\!f(n)$ and $B=Vect(n)$ or $Vect_0(n).$
Sometimes  the cup-product of $\psi\in C^k(A,B)$ and
$\phi\in C^l(A,B)$ lies in $C^{k+l}(B,B).$ Such fortunate situations
occur in calculating of $s_{k+l}$ for sufficiently large $k+l.$

\subsection{General case }
The difficult part of $k$-commutator calculations is related to the
problem of exact constructions of $k$-commutators. Calculations
for $k$-commutators are done mainly for the case $U={\bf
C}[x_1,\ldots,x_n].$ This allows us to use the escort method
effectively. After receiving exact formulas it is less
difficult to check that these formulas for $k$-commutators  are
valid for any $\mathcal D$-differential algebras. For example, for
the cases $U={\bf C}[x_1,\ldots,x_n]$ or $U=C^{\infty}(M)$ one can
use continuity properties of functions  and the fact that
polynomials are dense.

One can give a pure algebraic proof of these statements. It needs
some additional arguments based on the super variational calculus.
Our questions are equivalent to the following question. Let $D$ be
an odd derivation. Is it possible to find some $k,$ such that
$D^k$ is also derivation? For $k=2,$ it is well known.
Surprisingly, that for some $n$ one can find a special $k=k(n),$
as in case of $k$-commutators, such that $D^k$ will be
a derivation and $D^{k+1}=0.$ For example, $D^7=0$ if
$D=u_1\der_1+u_2\der_2,$ where $u_1$ and $u_2$ are odd variables
and $\der_1,\der_2$ are even derivations. In general,
$D^{n^2+2n-1}=0$ for $D=\sum_{i=1}^n u_i\der_i$ where $u_i$ are
odd elements and $\der_1,\ldots,\der_n$ are commuting even
derivations. It seems that this topic is interesting in itself, and we
are planning to discuss it in a separate paper.

We give here only some key ideas. Let $U$ be a ${\mathcal
D}$-differential algebra, e.g. the Laurent polynomial algebra
${\bf C}[x_1^{\pm 1},\ldots,x_n^{\pm 1}]$ or an algebra
of smooth functions on some $n$-dimensional manifold. We would
like to calculate $s_k$ on $Dif\!f(n)$ corresponding to this $U.$
Notice that all of these algebras contain ${\bf
C}[x_1,\ldots,x_n].$ Therefore, to have the identity $s_k=0$ on
$Vect(n),$ it is necessary to have the identity $s_k=0$ on $W_n.$ These
conditions in our cases are sufficient. In fact, let $U$ be
any $\D$-differential algebra and $X_l=u_l\der_{i_l}\in Vect(n)
\subseteq Dif\!f(n),$ where $u_l\in U$ and $i_l\in \{1,\ldots,n\}.$
If $i_l=s,$ we say that $X_l$ or $u_l$ has type $s.$
Notice that
$s_k(X_1,\ldots,X_k)$ is a sum of elements of the form
$\der^{\alpha^1}(u_1)\cdots
\der^{\alpha^k}(u_k)\der^{\alpha^{k+1}}$ where
$\alpha^1,\ldots,\alpha^k,\alpha^{k+1}\in {\bf Z}_+^n$ such that
$\alpha^1+\cdots+\alpha^{k+1}=k.$ Moreover, since $s_k$ is
skew-symmetric, if $\alpha^l=\alpha^s,$ then $u_l$ and $u_s$
should have different types. For example, monomials of the form
$u_{j_1}\ldots u_{j_l}$ should enter no more than $n$ times.

These kind of arguments shows that in obtaining the formula for
$s_k(X_1,\ldots,X_k)$ one uses only
\begin{itemize}
\item the Leibniz rule,
\item the polylinearity  of $s_k$,
\item the ${\mathcal D}$-invariance of $s_k$,
\item the $0$-grading property of $s_k.$
\end{itemize}
One can show that $s_k(u_1\der_{i_1},\ldots,u_k\der_{i_k})$ is a
sum of elements of the form
$\lambda_{(\alpha^1,u_1,\ldots,\alpha^k,u_k,\alpha^{k+1})}
\der^{\alpha^1}(u_1)\cdots
\der^{\alpha^k}(u_k)\der^{\alpha^{k+1}}$ where coefficients
$\lambda_{(\alpha^1,u_1,\ldots,\alpha^k,u_k,\alpha^{k+1})}$ do not
depend on the choice of $\mathcal D$-differential algebra $U.$
Therefore, we can get as a testing algebra a polynomials algebra.
This is the main way how to reduce calculations to the polynomials
case.

The following lemma \ref{askar} is stated in a general case,
although we prove it below essentially  in case of $U={\bf
C}[x_1,\ldots,x_n].$ One can check directly that statements of
this lemma are true independently of $U$ in each of the cases,
where we use these statements.

\begin{lm}\label{askar}  Let $s_k$ be a standard skew-symmetric polynomial
and $U$  an associative commutative algebra with a system of
commuting derivations ${\mathcal D}=\{\der_1,\ldots,\der_n\}.$ Let
$Dif\!f(n)=\langle u\der^{\alpha}: \alpha\in {\bf Z}^n_+, u\in U\rangle $ be a
space of differential operators on $U$ endowed with associative
multiplication $\cdot.$ Let $Vect(n)=\langle u\der_i: u\in U,
i=1,\ldots,n\rangle $ be the subspace of differential operators on $U,$ of
differential degree one. Let $W_n=\langle x^{\alpha}\der_i: \alpha\in
{\bf Z}^n+, i=1,\ldots,n\rangle $ be the Lie algebra corresponding to $U={\bf
C}[x_1,\ldots,x_n].$ Then

{\rm (i)} $s_{k,Dif\!f(n)}$ is ${\mathcal D}$-invariant.

{\rm (ii)} $s_{k,W_n}$ is $\mathcal E$-graded.

{\rm (iii)} $s_{k+l,Dif\!f(n)}=s_{k,Dif\!f(n)}\smile s_{l,Dif\!f(n)}.$

{\rm (iv)} $esc(s_{k+l,W_n}^{rsym.r}=s_{k,W_n}^{rsym.r}
\overset{\circ}{\smile} s_{l,Vect(n)},$ if $k+l\ge n^2+2n-2.$

{\rm (v)} $s_{k+l,S_1}^{rsym.r}=s_{k, S_1}^{rsym.r}\overset{\circ}
{\smile} s_{l, S_1},$ if $i+j\ge 5.$
\end{lm}

\section{Sufficient condition for a $\D$-invariant form with skew-symmetric
arguments to be zero}

Let $L$ be a graded algebra $L=\oplus_{i\ge -1}L_i$ with the commuting
system of derivations ${\mathcal D} =L_{-1}.$
Recall that  ${\mathcal M}=
\oplus_{s\ge -q}{\mathcal M}_{[s]}$ is a graded  $L$-module if
$$L_i \M_{[s]}\subseteq \M_{[i+s]},$$
for any $i\ge -1$ and $s\ge -q.$

Suppose $\mathcal M_1,\ldots,\mathcal M_k,\mathcal M_0$ are graded
$L$-modules. Say  $\psi\in C^k(\M_1,\ldots,\M_k;\M)$ is {\it graded} and
that $\psi$ has grade $r$ if
$$\psi(\M_{1,[s_1]},\ldots,\M_{k,[s_k]})\subseteq
\M_{0,[i_1+\cdots+i_k+r]},$$
for any $s_1,\ldots,s_k.$ Then
$$C^k(\M_1,\ldots,\M_k;\M_0)=
\oplus_s
C^k(\M_1,\ldots,\M_k;\M_0)_{[s]},$$
is also graded, where
$$C^k(\M_1,\ldots,\M_k;\M_0)_{[s]}=
<\psi\in C^k(\M_1,\ldots,\M_k;\M_0): |\psi|=s\rangle.$$

\begin{lm}\label{may14} Consider the linear programming problem
$$\left\{\begin{array}{l} \sum_{i=-1}^mx_i=r\\ 0<x_i\le l_i, \quad
i=-1,0,\ldots,m-1,\\
\\
f(x_{-1},x_{-1},\ldots,x_m):=\sum_{i=-1}^mix_i\rightarrow min.\\
\end{array}
\right.$$ Then $$\min f(x_{-1},x_{0},\ldots,x_m)=m
r-\sum_{i=-1}^{m-1}(m-i)x_i$$ and this value is attained for
$x_{-1}=l_{-1}, x_{0}=l_{0}, \ldots, x_{m-1}=l_{m-1}$ and
$x_m=r-\sum_{i=-1}^{m-1}l_i.$
\end{lm}

{\bf Proof.} Since $x_m=r-\sum_{i=-1}^{m-1}x_{i},$
$$f(x_{-1},x_{0},\ldots,x_m)=m\,r-\sum_{i=-1}^{m-1}(m-i)x_i.$$ Thus,
$$f(x_{-1},x_{0},\ldots,x_m)\le m\,r-
\sum_{i=-1}^{m-1}(m-i)l_i$$ and the inequality can be converted
to equality if $x_i=l_i, i=-1,0,\ldots,m-1,$ and
$x_m=r-\sum_{i=-1}^{m-1}l_i.$

\begin{thm}\label{may13} Let $A=\oplus_{i\ge -1} A_i$
and $M=\oplus_{i\ge q}M_i$ be  ${\mathcal D}$-graded
modules such that $A^{\mathcal D}=A_{-1}, M^{\mathcal D}=M_{q}.$
Suppose that $\psi\in T^k(A,M)$ is a $0$-graded polylinear
map and  skew-symmetric in $r$ arguments. Let $i_0$ be number
such that
$$\sum_{-1\le i\le i_0}\dim\,A_i\le r< \sum_{
-1\le i\le i_0+1}\dim\,A_i.$$
If
$$k+q< r(i_0+2)-\sum_{-1\le i\le
i_0}(i_0+1-i)\dim\,A_i$$
and $\psi$ is ${\mathcal D}$-invariant then $\psi=0.$
\end{thm}

{\bf Proof.} We prove that $esc(\psi)=0.$
Suppose that it is not true and $\psi\ne 0.$ Then there exist
homogeneous $a_1,\ldots,a_k\in A$ such that
$\psi(a_1,\ldots,a_k)\ne 0$ and $|a_1|+\cdots+|a_k|$ is minimal.
We have
$$\psi\in T^k(A,M)^{\mathcal D}\Rightarrow$$
$$\der_i\psi(a_1,\ldots,a_k)=\sum_{j=1}^k
\psi(a_1,\ldots,a_{j-1},\der_i(a_j),a_{j+1}, \ldots,a_k).$$
As
$|a_1|+\ldots+|a_k|$ is minimal with property
$\psi(a_1,\ldots,a_k)\ne 0$ and
$$|a_1|+\cdots+|a_k|>|a_1|+\cdots+|a_{j-1}|+|\der_i(a_j)|
+|a_{j+1}|+\cdots+|a_k|,$$ we obtain that
$$0\ne \psi(a_1,\ldots,a_k)\in M^{\mathcal D}=M_q.$$

Since $\psi$ is graded with degree 0, this means that we can
choose homogeneous elements $a_1,\ldots,a_k\in A$ such that
$$\psi(a_1,\ldots,a_k)\ne 0,\quad
|a_1|+\cdots+|a_k|=|\psi(a_1,\ldots,a_k)|=q.$$
As $\psi$ is
skew-symmetric in $r$ arguments, the set $\{a_1,\ldots,a_k\}$
should have at least $r$ linear independent elements. Denote them
by $a_{i_1},\ldots,a_{i_r}.$

Suppose that among $a_{i_1},\ldots,a_{i_r}$ there are $l_i$
elements of $A_i.$
Then \begin{equation} \label{arman10}
r=\sum_{i\ge -1}l_i \end{equation} and \begin{equation}
\label{arman11} l_i\le \dim\,A_i. \end{equation} Since $r\le
\sum_{i=-1}^{i_0+1}l_i,$ from (\ref{arman10}) it follows that
$$l_i=0, \quad i>i_0+1$$
and
\begin{equation}\label{arman12}
r=\sum_{i=-1}^{i_0+1}l_i.
\end{equation}
So, among elements $a_{i_1},\ldots,a_{i_r}$ there are $l_{-1}$
elements of degree $-1,$ $l_0$ elements of degree $0,$ etc,
$l_{i_0}$ elements of degree $i_0$ and finally
$r-\sum_{i=-1}^{i_0}l_i\ge
r-\sum_{i=-1}^{i_0}\dim\,A_i$ elements of degree $i_0+1.$ Since
$|a_i|\ge -1$ for any $i\ge -1,$ we obtain that
$$|\psi(a_1,\ldots,a_k)| =\sum_{i=1}^k|a_i|\ge (-1) (k-r)+
\sum_{s=-1}^{r} |a_{i_s}|\ge f(l_{-1},l_0,\ldots,l_{i_0+1}),$$
where $$ f(l_{-1},l_0,\ldots,l_{i_0+1})=
r-k+\sum_{i=-1}^{i_0+1}i\,l_i.$$ According to lemma~\ref{may14} and
our condition,
$$min\,f(l_{-1},l_0,\ldots,l_{i_0+1})=r-k+(i_0+1)r-
\sum_{i=-1}^{i_0+1} (i_0+1-i)\dim\,A_i> q. $$ Therefore,
$$|\psi(a_1,\ldots,a_k)|>q.$$ In particular,
$$\psi(a_1,\ldots,a_k)\not\in M^{\mathcal D},$$ which is a contradiction.

\begin{crl} \label{4-left commutative}
Let $\psi\in T^9(S_1,S_1)$ be a ${\mathcal D}$-invariant $0$-graded
form with $8$ skew-symmetric arguments. Then $\psi=0.$
\end{crl}

{\bf Proof.}
Recall that  $S_1$  denotes the  subalgebra of $W_2$ consisting
of derivations with divergence $0.$ Let $A=S_1.$ Then $\dim\,A_{-1}=2, \dim\,A_0=3,\dim\,A_1=4,$ since $A_{-1}=\langle\der_1,\der_2\rangle, $
$A_0=\langle x_1\der_2,x_1\der_1-x_2\der_2,x_2\der_1\rangle, $
$A_1=\langle x_2^2\der_1, x_1^2\der_1-2x_1x_2\der_2, x_2^2\der_2-
2x_1x_2\der_1, x_1^2\der_2\rangle.$
Thus,
$$\dim\,A_{-1}+\dim\,A_0=2+3\le r=8< \dim\,A_{-1}+\dim\,A_0+\dim\,A_1=2+3+4.$$
In other words, $i_0=0$ for $r=8.$ Furthermore, for $k=9, q=-1, r=8, i_0=0.$
We see that
$$k+q=8< 9=r(i_0+2)-\sum_{-1\le i\le i_0}(i_0+1-i)\dim\,A_i=8(0+2)-2 \cdot 2-1\cdot 3.$$
Therefore, all conditions of theorem \ref{may13} are fulfilled and $\psi=0$
for $A=S_1.$

\begin{crl} \label{5-left commutative} Let $\psi\in T^{11}(W_2, W_2)$
be a ${\mathcal D}$-invariant form with $10$ skew-symmetric arguments. Then
$\psi=0.$
\end{crl}

{\bf Proof.} Take $A=W_2$. Then $\dim\,A_{-1}=2, \dim\,A_{0}=4, \dim\,A_{1}=6,$
since
$A_{-1}=\langle\der_1,\der_2\rangle, $
$A_0=\langle x_i\der_j: i,j=1,2\rangle, $
$A_1=\langle x_i x_j\der_s: i,j,s=1,2, i\le j\rangle.$
For $r=10, k=11, q=-1$ we see that
$$\dim\,A_{-1}+\dim\,A_0=2+4\le r=10< \dim\,A_{-1}+\dim\,A_0+\dim\,A_1=2+4+6.$$
Hence $i_0=0,$ and
$$k+q=10< 12 =r(i_0+2)-\sum_{-1\le i\le i_0}(i_0+1-i)\dim\,A_i=
10(0+2)-2 \cdot 2-1\cdot 4.$$
Therefore, by theorem \ref{may13}, $\psi=0$ on $A=W_2.$

\begin{crl} \label{n^2+2n commutator}
$s_k=0$ is an identity on $W_n,$ if $k\ge n^2+2n.$
\end{crl}

{\bf Proof.}  Let $A=W_n.$ Then $\dim\,A_{-1}=n, \dim\,A_0=n^2,
\dim\,A_{1}=n^2(n+1)/2.$
We see that for $r=k\ge n^2+2n,$
$$\dim\,A_{-1}+\dim\,A_0=n+n^2\le r.$$
Therefore $i_0\ge 0.$ Hence, if $i_0=0$ then
$$k+q=r-1< r+r-2n-n^2=r(0+2)-2 \dim\,A_{-1}-\dim\,A_0.$$
If $i_0>0,$ $n>1,$ then $2\dim A_{-1}\le \dim A_1$ and
$$r(i_0+2)-\sum_{-1\le i\le i_0}(i_0+1-i)\dim\,A_i=$$
$$
r-2\dim A_{-1}+\sum_{1\le i\le i_0}i\dim A_i+
(i_0+1)(r-\sum_{-1\le i\le i_0}\dim\,A_{i})\ge $$
$$
r-2\dim A_{-1}+\sum_{1\le i\le i_0}i\dim A_i>r-1.$$
Thus,  we can use theorem \ref{may13}, if $n>1.$ If $n=1,$ it is easy to
check that $s_3=0,$ and, $s_k=0,$ for any $k>3.$

So, we have proved $\psi=0$ for $A=W_n.$

\begin{crl}\label{n^2+2n-3 left commutative}
Let $\psi\in T^{2n^2+4n-5}(W_n,W_n), n>1,$ be skew-symmetric in
$r\ge (3n^2+6n-5)/2$ arguments. Then $\psi=0.$
In particular, $(Vect(n),s_{n^2+2n-2})$ is $(n^2+2n-3)$-left commutative.
\end{crl}

{\bf Proof.} Let $A=W_n.$ For $q=-1,
k=2n^2+4n-5, r\ge (3n^2+6n-5)/2,$ it is easy to see that
$i_0\ge 0.$

Check that the case $i_0> 0$ is impossible.
If $n=2$ then we obtain a contradiction with the conditions
$$r\le k=11$$ and
$$\dim A_{-1}+\dim A_0+\dim A_1=12\le r.$$
Let $n>2.$ Then we will have
$$\dim A_{-1}+\dim A_0+\dim A_1=n+n^2+n^2(n+1)/2\le r <k=2n^2+4n-5,$$
and
$$n^3-n^2-6 n+5\le 0.$$
For $n\ge 3,$
$$n^3-n^2-6 n+5\ge 2n^2-6n+5>0,$$
and again obtain a contradiction.

So, $i_0=0.$ Then
$$k+q=2n^2+4n-6< 2 n^2+4n-5\le 2 r -2 \dim A_{-1}-\dim A_0.$$
Hence, the condition of theorem \ref{may13} is satisfied. Thus, $\psi=0$
for $A=W_n.$

Notice that $2(n^2+2n-3)\ge (3n^2+6n-5)/2$ if $n>1.$ Therefore,
the $(n^2+2n-3)$-left commutativity condition, as a condition
for a $\mathcal D$-invariant form
with $2(n^2+2n-3)$ skew-symmetric arguments, is an identity on $W_n.$

\section{Invariant $N$-operation on vector fields}

\begin{lm} \label{20fev}
The $sl_n$-module $\wedge^{n-2} (R(2\pi_1)\otimes R(\pi_{n-1})$
does not contain $R(2\pi_{n-1})$ as a submodule.
\end{lm}

{\bf Proof.} I am grateful to R.Howe for the following elegant proof of
this lemma.

One can argue that the full $(n-2)$ tensor
power of $R(2\pi_1)\otimes R(\pi_{n-1})$ does not contain
$R(2\pi_{n-1}).$ Indeed, the $(n-2)$ tensor power of this
tensor product is equal to the tensor product of the $(n-2)$
tensor powers of each factor.

The representation $R(2\pi_1)$
corresponds to the diagram with one row or length two. The
representation $R(\pi_{n-1})$ corresponds  to the diagram
with one column of length $n-1.$ So, the question then becomes,
does the Young diagram with $n-3$ columns of length $n,$ and
two columns of length $n-1,$ appear in the indicated tensor
product?

Since the diagram of $R(\pi_{n-1})$ has only one column, all of
the components of its $(n-2)$ tensor power will have at most
$(n-2)$ columns. Since the diagram of $R(2\pi_1)$ has only one
row, all the components of its $(n-2)$ tensor power will contain
at most $n-2$ rows. Now taking the tensor product of these two
representations, we can say that all components of the tensor
product will have diagrams which fit in an $\Gamma$-shaped region
with $(n-2)$ columns and $(n-2)$ rows. But the diagram of the
representation we are asking about does not fit in to this region,
so it cannot be a component.

\begin{crl}\label{28fe}
$s_{n^2+2n-2}$ has no quadratic differential part on
$Vect(n).$
\end{crl}

{\bf Proof.} Since, as $sl_n$-modules,
$$L_1\cong R(2\pi_1+\pi_{n-1})\oplus R(\pi_1)\cong
R(2\pi_1)\otimes R(\pi_{n-1}),$$
we obtain an isomorphism of $sl_n$-modules
$$\wedge^kL_1\cong \wedge^k(R(2\pi_1)\otimes R(\pi_{n-1})).$$

Consider the  homomorphism of $sl_n$-modules
$$\rho_{k,s}: \wedge^{k-n^2-n} L_1\rightarrow R(s\pi_{n-1}),$$
induced by
$$\rho(X_1\wedge \ldots\wedge X_{k-n^2-n})=$$
$$
pr_{s}(s_{k,Dif\!f(n)}
(\der_1,\ldots,\der_n,x_1\der_1,\ldots,x_n\der_1,\ldots,x_1\der_n,\ldots,x_n\der_n,X_1,\ldots,X_{k-n^2-n})),$$
where
$pr_{s}: Dif\!f(n)\rightarrow <\der^{\alpha}: |\alpha|=s\rangle\cong R(s\pi_{n-1})$
is the projection map.

Since
$$\wedge^nL_{-1}\otimes \wedge^{n^2}L_0\cong {\bf C},$$
it is clear that $\rho_{n^2+2n-2,2}$ should give a homomorphism of
$\wedge^{n-2}(R(\pi_1)\otimes R(\pi_1+\pi_{n-1}))$ to
$R(2\pi_{n-1}).$ By lemma \ref{20fev} this homomorphism is
trivial. Thus, $s_{n^2+2n-2}(X_1,\ldots,X_{n^2+2n-2})\in Vect(n)$
for any $X_1,\ldots,X_{n^2+2n-2}\in Vect(n).$

\begin{lm} \label{28fevruary} Let $U$ be any ${\mathcal D}$-differential
algebra. For any $X_1,\ldots,X_{k}\in Vect(n),$
$s_k(X_1,\ldots,X_k)\in Vect(n),$
if $k=n^2+2n-2$ and $k=n^2+2n-1.$
\end{lm}

{\bf Proof.} For $k=n^2+2n-2$ this follows from corollary \ref{28fe}.
For $k=n^2+2n-1$ we see that $esc(s_k)$ has support
$\wedge^{n}L_{-1}\otimes \wedge^{n^2-1}\otimes \wedge^{n-1}L_1$
and
$s_k(\der_1,\ldots,\der_n,a_1,\ldots,a_{n^-1},X_1,\ldots,X_{n-1}),$
where $a_1,\ldots,a_{n^2-1}\in L_0,$ $X_1,\ldots,X_{n-1}\in L_1,$ never
gives quadratic terms, as
$$|s_k(\der_1,\ldots,\der_n,a_1,\ldots,a_{n^-1},X_1,\ldots,X_{n-1})|=-1.$$

\begin{lm}\label{fev28}
$s_k=0$ is an identity on $Vect(n)$ for $k\ge n^2+2n.$
\end{lm}

{\bf Proof.} This is corollary \ref{n^2+2n commutator}. However, we give here
another proof.
 $s_k$ is skew-symmetric in $k\ge n^2+2n$ arguments and
graded. Since $\dim\,L_{-1}=n, \dim\,L_0=n^2, \dim\,L_1=n^2(n+1)/2\ge n$
and
$$|s_k(\der_1,\ldots,\der_n,x_1\der_1,\ldots,x_n\der_1,\ldots,x_1\der_n,
\ldots,x_n\der_n,X_1,\ldots,X_{k}))|\ge $$
$$ (k-n^2-n)\cdot 1+n^2\cdot 0+n\cdot(-1)\ge 0,$$
we see that $esc(s_k)=0,$ if $k\ge n^2+2n.$

{\bf Remark.} $n^2+2n$ here is not minimal. One can prove
that $s_{n^2+2n-1}=$ is an identity on $Vect(n)$ and that
$s_{n^2+2n-2}=0$ is an identity on $Vect_0(n).$

\section{Quadratic differential parts for $k$-commutators
in two variables}

Let $Diff_s(n)$ be a subspace of $Dif\!f(n)$  that consists of
differential operators of differential order $s,$ i.e.,
$Diff_s(n)=\langle u\der^{\alpha}: u\in U, \alpha\in {\bf Z}^n_+,
|\alpha|=\sum_{i}\alpha_i=s\rangle .$ Let $pr_s: Dif\!f(n)\rightarrow
Diff_s(n)$ be the projection map.

\begin{lm} For any $X_1,\ldots,X_k\in Vect(2),$
$$pr_l (s_k(X_1,\ldots,X_k))=0,$$
if $l>2.$
\end{lm}

{\bf Proof.} If $k>6$ then by lemma \ref{fev28} and
\ref{arman}, $s_k=0$ is an identity.
If $k=6,$ then $s_{k,Vect(2)}$ has only a linear part.
If $k\le 5$ then $s_k$ can be decomposed into a cup-product of $s_2$ and
$s_{k-2}.$ We know that $s_2$ can only give differential operators of first order. So, $s_3=s_2\smile s_1$ and $s_4=s_2\smile s_2$ can give differential
operators at most second order. As far as $s_5=s_3\smile s_2,$
the following reasoning  shows that the differential operators of third order
can not be represented as $s_5(X_1,\ldots,X_5)$ for any
$X_1,\ldots,X_5\in Vect(2).$ For $L=W_2,$  support for an escort map of $pr_ls_5$
with a maximal $l$ should contain $\{\der_1,\der_2,x_2\der_1,a,b,c\}$,
where $a,b,c\in L_0.$ Easy calculations then show that $l\le 2$ if
$k=5.$

{\bf Remark.} One can prove that if $l>n$ then
$pr_l(s_k(X_1,\ldots,X_k))=0$
for any $k$ and $X_1,\ldots,X_k\in Vect(n).$

\begin{lm}\label{quadratic part for s3}
$$pr_2(s_3(X_1,X_2,X_3))=$$
$$-\left|\begin{array}{ccc}
(x_1)X_1&(x_1)X_2&(x_1)X_3\\
(x_2)X_1&(x_2)X_2&(x_2)X_3\\
\der_2((x_1)X_1)&\der_2((x_1)X_2)&\der_2((x_1)X_3)\\
\end{array}\right|\der_1^2
$$
$$
+\left|\begin{array}{ccc}
(x_1)X_1&(x_1)X_2&(x_1)X_3\\
(x_2)X_1&(x_2)X_2&(x_2)X_3\\
\der_1((x_1)X_1) &\der_1((x_1)X_2)
&\der_1((x_1)X_3)\\
\end{array}\right|\der_1\der_2
$$
$$
-\left|\begin{array}{ccc}
(x_1)X_1&(x_1)X_2&(x_1)X_3\\
(x_2)X_1&(x_2)X_2&(x_2)X_3\\
\der_2((x_2)X_1) &\der_2((x_2)X_2)
&\der_2((x_2)X_3)\\
\end{array}\right|\der_1\der_2
$$
$$+\left|\begin{array}{ccc}
(x_1)X_1&(x_1)X_2&(x_1)X_3\\
(x_2)X_1&(x_2)X_2&(x_2)X_3\\
\der_1((x_2)X_1)&\der_1((x_2)X_2)&\der_1((x_2)X_3)\\
\end{array}\right|\der_2^2,
$$
for any $X_1,X_2,X_3\in Vect(2).$
\end{lm}

\begin{lm}\label{quadratic part for s4}
$$pr_2(s_4(X_1,\ldots,X_4))=$$
$$
-2\left|\begin{array}{cccc}
(x_1)X_1&(x_1)X_2&(x_1)X_3&(x_1)X_4\\
(x_2)X_1&(x_2)X_2&(x_2)X_3&(x_2)X_4\\
\der_1((x_1)X_1)&\der_1((x_1)X_2)&\der_1((x_1)X_3)&\der_1((x_1)X_4)\\
\der_2((x_1)X_1)&\der_2((x_1)X_2)&\der_2((x_1)X_3)&\der_2((x_1)X_4)\\
\end{array}\right|\der_1^2
$$
$$
-2\left|\begin{array}{cccc}
(x_1)X_1&(x_1)X_2&(x_1)X_3&(x_1)X_4\\
(x_2)X_1&(x_2)X_2&(x_2)X_3&(x_2)X_4\\
\der_1((x_1)X_1)&\der_1((x_1)X_2)&\der_1((x_1)X_3)&\der_1((x_1)X_4)\\
\der_2((x_2)X_1)&\der_2((x_2)X_2)&\der_2((x_2)X_3)&\der_2((x_2)X_4)\\
\end{array}\right|\der_1\der_2
$$
$$
-2\left|\begin{array}{cccc}
(x_1)X_1&(x_1)X_2&(x_1)X_3&(x_1)X_4\\
(x_2)X_1&(x_2)X_2&(x_2)X_3&(x_2)X_4\\
\der_1((x_2)X_1)&\der_1((x_2)X_2)&\der_1((x_2)X_3)&\der_1((x_2)X_4)\\
\der_2((x_1)X_1)&\der_2((x_1)X_2)&\der_2((x_1)X_3)&\der_2((x_1)X_4)\\
\end{array}\right|\der_1\der_2
$$
$$
-2\left|\begin{array}{cccc}
(x_1)X_1&(x_1)X_2&(x_1)X_3&(x_1)X_4\\
(x_2)X_1&(x_2)X_2&(x_2)X_3&(x_2)X_4\\
\der_1((x_2)X_1)&\der_1((x_2)X_2)&\der_1((x_2)X_3)&\der_1((x_2)X_4)\\
\der_2((x_2)X_1)&\der_2((x_2)X_2)&\der_2((x_2)X_3)&\der_2((x_2)X_4)\\
\end{array}\right|\der_2^2,
$$
for any $X_1,\ldots,X_4\in Vect(2).$
\end{lm}

\begin{lm} \label{quadratic part for s5}
$$pr_2(s_5(X_1,\ldots,X_5))=$$
$$
-\left|\begin{array}{ccccc}
(x_1)X_1&(x_1)X_2&(x_1)X_3&(x_1)X_4&(x_1)X_5\\
(x_2)X_1&(x_2)X_2&(x_2)X_3&(x_2)X_4&(x_2)X_5\\
\der_1((x_1)X_1)&\der_1((x_1)X_2)&\der_1((x_1)X_3)&\der_1((x_1)X_4)&\der_1((x_1)X_5)\\
\der_2((x_2)X_1)&\der_2((x_2)X_2)&\der_2((x_2)X_3)&\der_2((x_2)X_4)&
\der_2((x_2)X_5)\\
\der_2((x_1)X_1)&\der_2((x_1)X_2)&\der_2((x_1)X_3)&\der_2((x_1)X_4)&
\der_2((x_1)X_5)\\
\end{array}\right|\der_1^2
$$
$$
-\left|\begin{array}{ccccc}
(x_1)X_1&(x_1)X_2&(x_1)X_3&(x_1)X_4&(x_1)X_5\\
(x_2)X_1&(x_2)X_2&(x_2)X_3&(x_2)X_4&(x_2)X_5\\
Div\,X_1&Div\,X_2& Div\,X_3&Div\,X_4&Div\,X_5\\
\der_2((x_1)X_1)&\der_2((x_1)X_2)&\der_2((x_1)X_3)&\der_2((x_1)X_4)&
\der_2((x_1)X_5)\\
\der_1((x_2)X_1)&\der_1((x_2)X_2)&\der_1((x_2)X_3)&\der_1((x_2)X_4)&
\der_1((x_2)X_5)\\
\end{array}\right|\der_1\der_2
$$
$$
-\left|\begin{array}{ccccc}
(x_1)X_1&(x_1)X_2&(x_1)X_3&(x_1)X_4&(x_1)X_5\\
(x_2)X_1&(x_2)X_2&(x_2)X_3&(x_2)X_4&(x_2)X_5\\
\der_1((x_1)X_1)&\der_1((x_1)X_2)&\der_1((x_1)X_3)&\der_1((x_1)X_4)&
\der_1((x_1)X_5)\\
\der_2((x_2)X_1)&\der_2((x_2)X_2)&\der_2((x_2)X_3)&\der_2((x_2)X_4)&
\der_2((x_2)X_5)\\
\der_1((x_2)X_1)&\der_1((x_2)X_2)&\der_1((x_2)X_3)&\der_1((x_2)X_4)&
\der_1((x_2)X_5)\\
\end{array}\right|\der_2^2,
$$
for any $X_1,\ldots,X_5\in Vect(2).$
\end{lm}

To prove these statements one needs to  calculate their escorts.
A sufficient number of examples of similar calculations will be given below.

\section{Exact formula for $5$-commutator}

\begin{thm}\label{exact formula for s5}
Let $U$ be an associative commutative algebra with
two commuting derivations $\der_1$ and $\der_2.$ Then
$$s_5(D_{12}(u_1),D_{12}(u_2),D_{12}(u_3),D_{12}(u_4),D_{12}(u_5))=
-3 D_{12}([u_1,u_2,u_3,u_4,u_5]),$$
for any $u_1,\ldots,u_5\in U,$ where
$$[u_1,u_2,u_3,u_4,u_5]=
\left\vert\begin{array}{ccccc}
\der_1u_1&\der_1u_2&\der_1u_3&\der_1u_4&\der_1u_5\\
\der_2u_1&\der_2u_2&\der_2u_3&\der_2u_4&\der_2u_5\\
\der_1^2u_1&\der_1^2u_2&\der_1^2u_3&\der_1^2u_4&\der_1^2u_5\\
\der_1\der_2u_1&\der_1\der_2u_2&\der_1\der_2u_3&\der_1\der_2u_4&
\der_1\der_2u_5\\
\der_2^2u_1&\der_2^2u_2&\der_2^2u_3&\der_2^2u_4&\der_2^2u_5\\
\end{array}\right\vert$$
and $D_{12}(u)=\der_1(u)\der_2-\der_2(u)\der_1.$
\end{thm}

{\bf Proof.}
Let $L_i$ be graded components for $S_1=\langle x\in W_2: Div\,X=0\rangle $
and $a,b,c\in L_0, X\in L_1.$
Notice that
$$s_3^{rsym.r}(\der_i,a,X)=[\der_i,a]\circ X+[a,X]\circ\der_i+[X,\der_i]\circ a=$$
$$-a\circ\der_iX+X\circ \der_ia\in L_{0}.$$
Therefore, by lemma \ref{askar}
$$s_5(\der_1,\der_2,a,b,X)=$$

$$-s_3^{rsym.r}(\der_1,b,X)\circ \der_2(a)+s_3^{rsym.r}(\der_2,b,X)\circ
\der_1(a)$$
$$+s_3^{rsym.r}(\der_1,a,X)\circ \der_2(b)
-s_3^{rsym.r}(\der_2,a,X)\circ \der_1(b)=$$

$$+(b\circ\der_1X-X\circ\der_1b)\circ\der_2a-(b\circ \der_2X-X\circ\der_2b)
\circ \der_1a$$
$$+(X\circ\der_1a-a\circ\der_1X)\circ\der_2b-(X\circ \der_2a-a\circ\der_2X)
\circ \der_1b=$$

$$-(a\circ\der_1X-X\circ\der_1a)\circ\der_2b+
(a\circ\der_2X-X\circ\der_2a)\circ\der_1b$$
$$+(b\circ\der_1X-X\circ\der_1b)\circ\der_2a
-(b\circ\der_2X-X\circ\der_2b)\circ\der_1a=$$

$$-a\circ[\der_1X,\der_2b]+[X,\der_2b]\circ\der_1a
+a\circ[\der_2X,\der_1b]-[X,\der_1b]\circ\der_2a$$
$$+b\circ[\der_1X,\der_2a]-[X,\der_2a]\circ\der_1b
-b\circ[\der_2X,\der_1a]+[X,\der_1a]\circ\der_2b.$$
We see that non-zero components for $esc(s_5)$ are
$$s_5(\der_1,\der_2,x_2\der_1,x_1\der_1-x_2\der_2, x^2_1\der_2)=6\der_2,$$
$$s_5(\der_1,\der_2,x_2\der_1,x_1\der_1-x_2\der_2, x^2_1\der_1-2 x_1x_2\der_2)
=6\der_1,$$
$$s_5(\der_1,\der_2,x_2\der_1,x_1\der_2, x^2_1\der_1-2 x_1x_2\der_2)
=-6\der_2,$$
$$s_5(\der_1,\der_2,x_2\der_1,x_1\der_2, x^2_2\der_2-2 x_1x_2\der_1)
=-6\der_1,$$
$$s_5(\der_1,\der_2,x_1\der_1-x_2\der_2,x_1\der_2, x^2_2\der_2
-2 x_1x_2\der_1)=-6\der_2,$$
$$s_5(\der_1,\der_2,x_1\der_1-x_2\der_2,x_1\der_2, x^2_2\der_1)=-6\der_1.$$
It is easy to check, that
$$esc(s_5^{rsym.r})(D_{12}(u_1),\ldots,D_{12}(u_5))=-3 pr_{-1} D_{12}([u_1,\ldots,u_5]),$$
for any $u_1,\ldots,u_5\in {\bf C}[x_1,x_2],$ such that $|u_1|+\cdots+|u_5|=11.$

It remains to use (\ref{escort}) for $\mathcal D$-invariant form
$s_5^{rsym.r}$ and use lemma~\ref{quadratic part for s5}.

\section{Exact formula for $6$-commutator}

In this section we give a formula for calculating
the  $6$-commutator $s_6(X_1,\ldots,X_6)$
for $X_i=u_{i,1}\der_1+u_{i,2}\der_2,$
$i=1,\ldots,6.$
It can be presented as a sum of fourteen $6\times 6$ determinants
of the form
$$\lambda\left|\begin{array}{cccccc}
u_{1,1}&u_{2,1}&u_{3,1}&u_{4,1}&u_{5,1}&u_{6,1}\\
u_{1,2}&u_{2,2}&u_{3,2}&u_{4,2}&u_{5,2}&u_{6,2}\\
*&*&*&*&*&*\\
\end{array}\right|\der_i$$
where $\lambda =-1,2,-3$ and $i=1,2.$ Write $4\times 6$ parts of
such matrices  denoted by $*$ for $i=1.$  They are
$$-\left|\begin{array}{cccccc}
\der_2u_{1,1}&\der_2u_{2,1}&\der_2u_{3,1}&
\der_2u_{4,1}&\der_2u_{5,1}&\der_2u_{6,1}\\
\der_1u_{1,2}&\der_1u_{2,2}&\der_1u_{3,2}&
\der_1u_{4,2}&\der_1u_{5,2}&\der_1u_{6,2}\\
\der_2u_{1,2}&\der_2u_{2,2}&\der_2u_{3,2}&
\der_2u_{4,2}&\der_2u_{5,2}&\der_2u_{6,2}\\
\der_2^2u_{1,2}&\der_2^2u_{2,2}&\der_2^2u_{3,2}&
\der_2^2u_{4,2}&\der_2^2u_{5,2}&\der_2^2u_{6,2}\\
\end{array}\right|\der_1
$$

$$-\left|\begin{array}{cccccc}
\der_1u_{1,1}&\der_1u_{2,1}&\der_1u_{3,1}&
\der_1u_{4,1}&\der_1u_{5,1}&\der_1u_{6,1}\\
\der_2u_{1,1}&\der_2u_{2,1}&\der_2u_{3,1}&
\der_2u_{4,1}&\der_2u_{5,1}&\der_2u_{6,1}\\
\der_2u_{1,2}&\der_2u_{2,2}&\der_2u_{3,2}&
\der_2u_{4,2}&\der_2u_{5,2}&\der_2u_{6,2}\\
\der_1^2u_{1,1}&\der_1^2u_{2,1}&\der_1^2u_{3,1}&
\der_1^2u_{4,1}&\der_1^2u_{5,1}&\der_1^2u_{6,1}\\
\end{array}\right|\der_1
$$

$$-\left|\begin{array}{cccccc}
\der_1u_{1,1}&\der_1u_{2,1}&\der_1u_{3,1}&
\der_1u_{4,1}&\der_1u_{5,1}&\der_1u_{6,1}\\
\der_2u_{1,1}&\der_2u_{2,1}&\der_2u_{3,1}&
\der_2u_{4,1}&\der_2u_{5,1}&\der_2u_{6,1}\\
\der_1u_{1,2}&\der_1u_{2,2}&\der_1u_{3,2}&
\der_1u_{4,2}&\der_1u_{5,2}&\der_1u_{6,2}\\
\der_2^2u_{1,2}&\der_2^2u_{2,2}&\der_2^2u_{3,2}&
\der_2^2u_{4,2}&\der_2^2u_{5,2}&\der_2^2u_{6,2}\\
\end{array}\right|\der_1
$$

$$+2\left|\begin{array}{cccccc}
\der_2u_{1,1}&\der_2u_{2,1}&\der_2u_{3,1}&
\der_2u_{4,1}&\der_2u_{5,1}&\der_2u_{6,1}\\
\der_1u_{1,2}&\der_1u_{2,2}&\der_1u_{3,2}&
\der_1u_{4,2}&\der_1u_{5,2}&\der_1u_{6,2}\\
\der_2u_{1,2}&\der_2u_{2,2}&\der_2u_{3,2}&
\der_2u_{4,2}&\der_2u_{5,2}&\der_2u_{6,2}\\
\der_{12}u_{1,1}&\der_{12}u_{2,1}&\der_{12}u_{3,1}&
\der_{12}u_{4,1}&\der_{12}u_{5,1}&\der_{12}u_{6,1}\\
\end{array}\right|\der_1
$$

$$+2\left|\begin{array}{cccccc}
\der_1u_{1,1}&\der_1u_{2,1}&\der_1u_{3,1}&
\der_1u_{4,1}&\der_1u_{5,1}&\der_1u_{6,1}\\
\der_2u_{1,1}&\der_2u_{2,1}&\der_2u_{3,1}&
\der_2u_{4,1}&\der_2u_{5,1}&\der_2u_{6,1}\\
\der_1u_{1,2}&\der_1u_{2,2}&\der_1u_{3,2}&
\der_1u_{4,2}&\der_1u_{5,2}&\der_1u_{6,2}\\
\der_{12}u_{1,1}&\der_{12}u_{2,1}&\der_{12}u_{3,1}&
\der_{12}u_{4,1}&\der_{12}u_{5,1}&\der_{12}u_{6,1}\\
\end{array}\right|\der_1
$$

$$+2\left|\begin{array}{cccccc}
\der_1u_{1,1}&\der_1u_{2,1}&\der_1u_{3,1}&
\der_1u_{4,1}&\der_1u_{5,1}&\der_1u_{6,1}\\
\der_2u_{1,1}&\der_2u_{2,1}&\der_2u_{3,1}&
\der_2u_{4,1}&\der_2u_{5,1}&\der_2u_{6,1}\\
\der_2u_{1,2}&\der_2u_{2,2}&\der_2u_{3,2}&
\der_2u_{4,2}&\der_2u_{5,2}&\der_2u_{6,2}\\
\der_{12}u_{1,2}&\der_{12}u_{2,2}&\der_{12}u_{3,2}&
\der_{12}u_{4,2}&\der_{12}u_{5,2}&\der_{12}u_{6,2}\\
\end{array}\right|\der_1
$$

$$-3\left|\begin{array}{cccccc}
\der_1u_{1,1}&\der_1u_{2,1}&\der_1u_{3,1}&
\der_1u_{4,1}&\der_1u_{5,1}&\der_1u_{6,1}\\
\der_1u_{1,2}&\der_1u_{2,2}&\der_1u_{3,2}&
\der_1u_{4,2}&\der_1u_{5,2}&\der_1u_{6,2}\\
\der_2u_{1,2}&\der_2u_{2,2}&\der_2u_{3,2}&
\der_2u_{4,2}&\der_2u_{5,2}&\der_2u_{6,2}\\
\der_2^2u_{1,1}&\der_2^2u_{2,1}&\der_2^2u_{3,1}&
\der_2^2u_{4,1}&\der_2^2u_{5,1}&\der_2^2u_{6,1}\\
\end{array}\right|\der_1
$$
Here we use the notation $\der_{12}=\der_1\der_2.$
Corresponding matrices for the $\der_2$ parts can be obtained from
these matrices by changing  $1$ to $2$ and $2$ to $1.$

\begin{thm}\label{s6}
$6$-commutator $s_6$ on $Vect(2)$ is be given by the formula given above.
\end{thm}

{\bf Proof.}
Let $X_1=\der_1,X_2=\der_2,$
$a_1=x_1\der_1, a_2=x_2\der_1,a_3=x_1\der_2,a_4=x_2\der_2.$
Let $V$ be the set of $6$-tuples of the form
$(X_1,X_2,a_1,\ldots\hat{a_i},\ldots,a_4,X_6),$ where $i=1,2,3,4$ and
$X_6$ runs over the basic elements of $W_2$ with order $|X_6|=1.$

We see that
$$supp(s_6)\subseteq V.$$

One calculates that
$$s_6(\der_1,\der_2,x_i\der_i,x_2\der_1,x_1\der_2,x_1^2\der_1)=-2 \der_2,$$
$$s_6(\der_1,\der_2,x_i\der_i,x_2\der_1,x_1\der_2,x_1 x_2\der_1)=2 \der_1,$$
$$s_6(\der_1,\der_2,x_i\der_i,x_2\der_1,x_1\der_2,x_1 x_2\der_2)=2\der_2,$$
$$s_6(\der_1,\der_2,x_i\der_i,x_2\der_1,x_1\der_2,x_2^2\der_2)=-2\der_1,$$
for $i=1,2$ and
$$s_6(\der_1,\der_2,x_1\der_1,x_2\der_1,x_2\der_2,x_1^2\der_1)=-2\der_1,$$
$$s_6(\der_1,\der_2,x_1\der_1,x_2\der_1,x_2\der_2,x_1 x_2 \der_2)=2\der_1,$$
$$s_6(\der_1,\der_2,x_1\der_1,x_2\der_1,x_2\der_2,x_1^2\der_2)=-6\der_2,$$

$$s_6(\der_1,\der_2,x_1\der_1,x_1\der_2,x_2\der_2,x_1 x_2\der_1)=2 \der_2,$$
$$s_6(\der_1,\der_2,x_1\der_1,x_1\der_2,x_2\der_2,x_2^2\der_1)=-6 \der_1,$$
$$s_6(\der_1,\der_2,x_1\der_1,x_1\der_2,x_2\der_2,x_2^2\der_2)=-2 \der_2,$$
For other $(X_1,\ldots,X_6)\in V,$
$$s_6(X_1,\ldots,X_6)=0.$$

Calculations here are not difficult, but tedious. We perform them  in
one example.

Take $$X_1=\der_1, X_2=\der_2, X_3=x_2\der_1, X_4=x_1\der_1-x_2\der_2,
X_5=x_1\der_2, X_6=x_1^2\der_1.$$
Then
$$s_6''(X_1,\ldots,X_5,x_1^2\der_1)=$$
$$2 (3 x_1^2 \der_1+2 x_1x_2\der_2)\circ \der_1\der_2
+2 x_1x_2\der_1 \circ \der_2^2=4\der_2.$$

We see  that
$$s_3^{rsym}(X_1,X_3,X_6)=s_3(\der_1,x_2\der_1,x_1^2\der_1)=
2 x_2\der_1\circ x_1\der_1=0,$$
and
$$s_3^{rsym}(X_1,X_3,X_6)\circ s_3^{rsym}(X_2,X_4,X_5)=0.$$
Furthermore,
$$s_3^{rsym}(X_1,X_4,X_6)=s_3^{rsym}(\der_1,x_1\der_1-x_2\der_2,x_1^2\der_1)=$$
$$2 (x_1\der_1-x_2\der_2)\circ x_1\der_1-x_1^2\der_1\circ \der_1=
2 x_1\der_1-2x_1\der_1=0,$$
and
$$s_3^{rsym}(X_1,X_4,X_6)\circ s_3^{rsym}(X_2,X_3,X_5)=0.$$
At last,
$$s_3^{rsym}(X_1,X_5,X_6)=s_3^{rsym}(\der_1,x_1\der_2,x_1^2\der_1)=$$
$$2 x_1\der_2\circ x_1\der_1-x_1^2\der_1\circ \der_2=2 x_1\der_2,$$

$$s_3^{rsym}(X_2,X_3,X_4)=s_3^{rsym}(\der_2,x_2\der_1,x_1\der_1-x_2\der_2)=$$
$$-x_2\der_1\circ \der_2-(x_1\der_1-x_2\der_2)\circ \der_1=-2\der_1,$$
and
$$s_3^{rsym}(X_1,X_5,X_6)\circ s_3^{rsym}(X_2,X_3,X_4)=-4\der_2.$$

Similarly,
$$s_3^{rsym}(X_2,X_3,X_6)=s_3^{rsym}(\der_2,x_2\der_1,x_1^2\der_1)=
-x_1^2\der_1\circ \der_1=-2x_1\der_1,$$

$$s_3^{rsym}(X_1,X_4,X_5)=s_3^{rsym}(\der_1,x_1\der_1-x_2\der_2,x_1\der_2)=$$
$$(x_1\der_1-x_2\der_2)\circ \der_2-x_1\der_2\circ \der_1= -2\der_2,$$
and
$$s_3^{rsym}(X_2,X_3,X_6)\circ s_3^{rsym}(X_1,X_4,X_5)=0.$$
We have:
$$s_3^{rsym}(X_2,X_4,X_6)=s_3^{rsym}(\der_2,x_1\der_1-x_2\der_2,x_1^2\der_1)=$$
$$x_1^2\der_1\circ \der_2=0,$$
and
$$s_3^{rsym}(X_2,X_4,X_6)\circ s_3^{rsym}(X_1,X_3,X_5)=0.$$
Finally,
$$s_3^{rsym}(X_2,X_5,X_6)=s_3^{rsym}(\der_2,x_1\der_2,x_1^2\der_1)=0,$$
and
$$s_3^{rsym}(X_2,X_5,X_6)\circ s_3^{rsym}(X_1,X_3,X_4)=0.$$
Thus,
$$s_6'(X_1,\ldots,X_5,x_1^2\der_1)=
s_3^{rsym}(X_1,X_5,X_6)\circ s_3^{rsym}(X_2,X_3,X_4)=-4\der_2.$$
Hence
$$s_6(X_1,\ldots,X_5,x_1^2\der_1)=
s_6'(X_1,\ldots,X_5,x_1^2\der_1)+
s_6''(X_1,\ldots,X_5,x_1^2\der_1)=
$$
$$-4\der_2+4\der_2=0.$$

So, we have constructed $esc(s_{6,W_2}^{rsym.r}).$
A reconstruction of $s_{6,W_2}^{rsym.r}$ by its escort
(formula (\ref{escort})) gives us the formula for $s_6.$
By lemma \ref{28fevruary}, $s_6=s_6^{rsym.r}$ on $Vect(2).$

\section{5-commutator of adjoint derivations}

\begin{lm} \label{s5 for adjoint} Let $U$ be any
$\{\der_1,\der_2\}$-differential algebra and $Vect_0(2)$
be the subspace of vector fields without divergence of $Vect(2).$
Then
$$ad\,s_5(X_1,\ldots,X_5)= s_5(ad\,X_1,\ldots,ad\,X_5),$$
for any $X_1,\ldots,X_5\in Vect_0(2).$
\end{lm}

{\bf Proof.}
Consider a multilinear polynomial $f$ with 6 variables defined by
$$f(t_0,t_1,\ldots,t_5)=(t_0)s_5^{ad}(t_1,\ldots,t_5)=
\sum_{\sigma\in Sym_5}sign\,\sigma\,
[\cdots[[t_0,t_{\sigma(1)}],t_{\sigma(2)}],\cdots].$$

We see that $f$ is polylinear and skew-symmetric in all variables
except the first one. Important properties for us are: $f_{Vect_0(2)}$ is
$\mathcal D$-invariant and $\mathcal E$-graded. Therefore, $f$ can be
uniquely restaured from its escort. We see that
$$supp(f)=$$
$$L_{-1}\otimes \wedge^2 L_{-1}\otimes \wedge^2L_0 \otimes L_2 \;\; \oplus\;\;  L_{-1}\otimes \wedge^2L_{-1}\otimes L_0\otimes \wedge^2L_1$$
$$\oplus\;\;L_{-1}\otimes L_{-1}\otimes \wedge^3L_0\otimes L_1$$
$$\oplus \;\;L_{0}\otimes \wedge^2L_{-1}\otimes \wedge^2L_0\otimes L_1$$
$$\oplus\;\;L_1\otimes \wedge^2L_{-1}\otimes \wedge^3L_0.$$
Here by $L_i$ we denote the graded components for
$S_1=\langle X\in W_2: Div\,X=0\rangle.$

Let $(a,b,c)=(x_2\der_1,x_1\der_1-x_2\der_2,x_1\der_2).$ Then
$$F:=s_5^{ad}(\der_1,\der_2,a,b,c)]\in End\,W_n,$$
is defined by
\begin{equation}\label{01march}
(u_1\der_1+u_2\der_2)F=
6(\der_1\der_2(u_1)+6\der_2^2(u_2))\der_1-6(\der_1^2u_1+\der_1\der_2u_2)\der_2.
\end{equation}
In other words,
$$(u_1\der_1+u_2\der_2)F=
-6D_{12}(\der_1(u_1)+\der_2(u_2)).$$

Set $s_k^{ad}(X_1,\ldots,X_k)=s_k(ad\,X_1,\ldots,ad\,X_k)),$ where
$ad: L\rightarrow End\,L.$

Let us prove (\ref{01march}).  By lemma \ref{askar} we have
$$F=F_1+F_2+F_3,$$
where
$$F_1=s_3^{ad}(\der_1,a,b){\cdot} ad[\der_2,c]+
s_3^{ad}(\der_1,b,c){\cdot} ad[\der_2,a]+
s_3^{ad}(\der_1,c,a){\cdot} ad[\der_2,b],$$
$$F_2=-s_3^{ad}(\der_2,a,b)\cdot ad[\der_1,c]-
s_3^{ad}(\der_2,b,c)\cdot ad[\der_1,a]-
s_3^{ad}(\der_2,c,a)\cdot ad[\der_1,b],$$
$$F_3=s_3^{ad}(\der_1,\der_2,a)\cdot ad\,[b,c]+
s_3^{ad}(\der_1,\der_2,b)\cdot ad[c,a]+
s_3^{ad}(\der_1,\der_2,c)\cdot ad[a,b].$$

Further,
$$F_1=-s_3^{ad}(\der_1,b,c)\cdot \der_1+
s_3^{ad}(\der_1,c,a)\cdot \der_2,$$
$$F_2=s_3^{ad}(\der_2,a,b)\cdot \der_2+
s_3^{ad}(\der_2,c,a)\cdot \der_1,$$
$$F_3=-2\,s_3^{ad}(\der_1,\der_2,a)\cdot  ad\,c-
s_3^{ad}(\der_1,\der_2,b)\cdot ad\,b
-2\,s_3^{ad}(\der_1,\der_2,c)\cdot ad\,a.$$

It is easy to see that
$$F_3=$$

$$2(\der_1\cdot \der_2a)\cdot ad\,c
-2(\der_2\cdot\der_1a)\cdot ad\,c$$
$$+(\der_1\cdot \der_2b)\cdot ad\,b
-(\der_2\cdot\der_1b)\cdot ad\,b$$
$$+2(\der_1\cdot \der_2c)\cdot ad\,a
-2(\der_2\cdot\der_1c)\cdot ad\,a=$$

$$2(\der_1\cdot \der_1)\cdot ad\,c-(\der_1\cdot \der_2)\cdot ad\,b
-(\der_2\cdot\der_1)\cdot ad\,b-2(\der_2\cdot\der_2)\cdot ad\,a=$$

$$2\der_1^2\cdot ad\,c
-2(\der_1 \der_2)\cdot ad\,b
-2\der_2^2\cdot ad\,a.$$
Note that
$$((u_1\der_1+u_2\der_2)\der_1^2) ad\,c=$$
$$[\der_1^2(u_1)\der_1+\der_1^2(u_2)\der_2,x_1\der_2]=
-\der_1^2(u_1)\der_2,$$

$$(u_1\der_1+u_2\der_2)(\der_1\der_2\cdot ad\,b)=$$
$$[\der_1\der_2(u_1)\der_1+\der_1\der_2(u_2)\der_2,
x_1\der_1-x_2\der_2]=0,$$

$$((u_1\der_1+u_2\der_2)\der_2^2) ad\,a=$$
$$[\der_2^2(u_1)\der_1+\der_2^2(u_2)\der_2,x_2\der_1]=
-\der_2^2(u_2)\der_1.$$
Therefore,
$$(u_1\der_1+u_2\der_2)F_3=$$
$$-2\der_1^2(u_1)\der_2+2\der_2^2(u_2)\der_1.$$

Similarly,
$$(u_1\der_1+u_2\der_2)(F_1+F_2)=$$
$$4\der_1\der_2(u_1)\der_1-4\der_1^2(u_1)\der_2+4\der_2^2(u_2)\der_1
-4\der_1\der_2(u_2)\der_2.$$

From these expressions of $F_1,$ $F_2$ and $F_3$ we obtain
(\ref{01march}).

So, by (\ref{01march}), $F=0$ on $Vect_0(2).$ In particular
$$esc(f)(X,\der_1,\der_2,x_2\der_1,x_1\der_1-x_2\der_2,
x_1\der_2)=0,$$
for any $X\in L_1.$

Similar calculations show that non-zero components for $esc(f)$
are
$$f(x_2\der_1,\der_1,\der_2,x_2\der_1,
x_1\der_1-x_2\der_2, D_{12}(x_1^3))=18 \der_1,$$
$$f(x_2\der_1,\der_1,\der_2,x_2\der_1,
x_1\der_2, D_{12}(x_1^2 x_2))=6 \der_1,$$
$$f(x_2\der_1,\der_1,\der_2,x_1\der_1-x_2\der_2,
x_1\der_2, D_{12}(x_1x_2^2))=-6 \der_1,$$

$$f(x_1\der_1-x_2\der_2,\der_1,\der_2,x_2\der_1,
x_1\der_1-x_2\der_2, D_{12}(x_1^3))=-18 \der_2,$$
$$f(x_1\der_1-x_2\der_2,\der_1,\der_2,x_2\der_1,
x_1\der_1-x_2\der_2, D_{12}(x_1^2 x_2))=-6 \der_1,$$
$$f(x_1\der_1-x_2\der_2,\der_1,\der_2,x_1\der_1-x_2\der_2,x_1\der_2, D_{12}(x_1 x_2^2))=6 \der_2,$$
$$f(x_1\der_1-x_2\der_2,\der_1,\der_2,x_1\der_1-x_2\der_2,x_1\der_2, D_{12}(x_2^3))=18\der_1,$$

$$f(x_1\der_2,\der_1,\der_2,x_2\der_1,
x_1\der_2, D_{12}(x_1 x_2^2))=-6 \der_2,$$
$$f(x_1\der_2,\der_1,\der_2,x_2\der_1,
x_1\der_1-x_2\der_2, D_{12}(x_1^2 x_2))=-6 \der_2,$$
$$f(x_1\der_2,\der_1,\der_2,
x_1\der_1-x_2\der_2,x_1\der_2, D_{12}(x_2^3))= 18 \der_2,$$

\medskip

\noindent and

\medskip

$$f(\der_1,\der_1,\der_2,x_2\der_1,x_1\der_1-x_2\der_2, D_{12}(x_1^4))=
-72\der_2,$$
$$f(\der_1,\der_1,\der_2,x_2\der_1,x_1\der_1-x_2\der_2,D_{12}(x^3_1 x_2))=18\der_1,$$
$$f(\der_1,\der_1,\der_2,x_2\der_1,x_1\der_2,D_{12}(x_1^3 x_2)=
-18 \der_2,$$
$$f(\der_1,\der_1,\der_2,x_2\der_1,x_1\der_2,D_{12}(x_1^2 x_2^2)=
12 \der_1,$$
$$f(\der_1,\der_1,\der_2,x_1\der_1-x_2\der_2,x_1\der_2, D_{12}(x_1^2 x_2^2))=12 \der_2,$$
$$f(\der_1,\der_1,\der_2,x_1\der_1-x_2\der_2,x_1\der_2, D_{12}(x_1 x_2^3))=-18 \der_1,$$

$$f(\der_1,\der_1,x_2\der_1,x_1\der_1-x_2\der_2, x_1\der_2,
D_{12}(x_1^3))=-18 \der_2,$$
$$f(\der_1,\der_1,x_2\der_1,x_1\der_1-x_2\der_2, x_1\der_2,
D_{12}(x_1^2x_2))=6 \der_1,$$
$$f(\der_1,\der_2,x_2\der_1,x_1\der_1-x_2\der_2, x_1\der_2,
D_{12}(x_1^2x_2))=-6 \der_2,$$
$$f(\der_1,\der_2,x_2\der_1,x_1\der_1-x_2\der_2, x_1\der_2,
D_{12}(x_1 x_2^2))=6\der_1,$$

$$f(\der_1,\der_1,\der_2,x_2\der_1,D_{12}(x_1^3), D_{12}(x_1^2x_2))=
-72\der_2,$$
$$f(\der_1,\der_1,\der_2,x_2\der_1,D_{12}(x_1^3), D_{12}(x_1x_2^2))=
36\der_1,$$
$$f(\der_1,\der_1,\der_2,x_1\der_1-x_2 \der_2,D_{12}(x_1^3),
D_{12}(x_1x_2^2))=72\der_2,$$
$$f(\der_1,\der_1,\der_2,x_1\der_1-x_2\der_2,D_{12}(x_1^3), D_{12}(x_2^3))=-108\der_1,$$
$$f(\der_1,\der_1,\der_2,x_1\der_2,D_{12}(x_1^2 x_2), D_{12}(x_2^3))=
-36\der_1,$$
$$f(\der_1,\der_1,\der_2,x_1\der_2,D_{12}(x_1^2 x_2), D_{12}(x_1x_2^2))=24\der_2.$$
Other components of the form $f(\der_1, \der_1,\der_2,a,b,X)$
$f(\der_1,\der_i,a,b,c,Y),$ and $f(\der_1,\der_1,\der_2,a,Y,Z)$
are equal to $0,$ where $a,b,c\in L_0,$  $X\in L_2, Y,Z\in L_1$ and
$i=1,2.$ To find
$f(\der_2,\der_1,\der_2,a,b,X),$ $f(\der_1,\der_i,a,b,c,Y)$
and $f(\der_2,\der_1,\der_2,a,Y,Z)$
one should use  the involutive
automorphism of $W_n$ induced by changing of variables
$(x_1,x_2)\mapsto (x_2,x_1).$

Using lemma \ref{exact formula for s5} we see that
$$esc(f)=esc(g),$$
where $\mathcal D$-invariant map $g: S_1\otimes \wedge^5S_1\rightarrow S_1$ is given by
$$g(X_1,X_1,\ldots,X_5)=[X_0,s_5^{rsym.r}(X_1,\ldots,X_5)].$$

It remains to use (\ref{escort}), for $\mathcal D$-invariant forms $f$ and $g$
to obtain that
$s_5^{ad}=ad\,s_5$ for $Vect_0(2).$

\section{$s_6=0$ is an identity on $Vect_0(2)$}

\begin{lm}\label{bir}
$s_6=0$ is an identity on $Vect_0(2).$
\end{lm}

{\bf Proof.} Set
$$X_1=\der_1, X_2=\der_2, X_3=x_2\der_1, X_4=x_1\der_1-x_2\der_2,
X_5=x_1\der_2,$$
$$V=\{(X_1,X_2,\ldots,X_6): |X_6|=1, X_6\in S_1\}.$$
Notice that $supp(s_6) \subseteq V.$ We need to check that
$s_6(X_1,\ldots,X_6)=0,$ for all $(X_1,\ldots,X_6)\in V.$ By lemma
\ref{28fevruary},
$$s_6^{rsym.r}=s_6.$$
By lemma \ref{askar},
$$s_6=s_3\smile s_3.$$

Let $(X_1,\ldots,X_6)\in V$ and $F=s_6(X_1,\ldots,X_6).$

We see that, $s_6(X_1,\ldots,X_6)$ is the alternating sum of
elements of the form
$s_3(X_{\sigma(1)},X_{\sigma(2)},X_{\sigma(3)})\cdot
s_3(X_{\sigma(4)},X_{\sigma(5)},X_{\sigma(6)}),$ where $\sigma\in
Sym_{3,3}$ are shuffle permutations, i.e.,
$\sigma(1)<\sigma(2)<\sigma(3), \sigma(4)<\sigma(5)<\sigma(6).$
Moreover,
$$s_6(X_1,\ldots,X_6)=$$
$$\sum_{\sigma\in Sym_{3,3}, \sigma(1)<\sigma(4)}
sign\,\sigma\,
[s_3(X_{\sigma(1)},X_{\sigma(2)},X_{\sigma(3)}),
s_3(X_{\sigma(4)},X_{\sigma(5)},X_{\sigma(6)})].$$

Since $|s_6(X_1,\ldots,X_6)|=-1,$
$$s_6(X_1,\ldots,X_6)\in <\der_i\rangle, $$
for some $i=1,2.$ Therefore, in calculating $F=s_6(X_1,\ldots,X_6)$
we can make summation only in
$\sigma\in Sym_{3,3}$ such that
$$s_3(X_{\sigma(1)}, X_{\sigma(2)},X_{\sigma(3)})\in <u\der_i: |u|=1,2\rangle, $$
$$s_3(X_{\sigma(4)}, X_{\sigma(5)},X_{\sigma(6)})\in <\der^{\alpha}:
|\alpha|=1,2\rangle.$$

Since
$$s_3(X_1,X_2,X)=s_3(\der_1,\der_2,X)=[\der_2, X]\cdot \der_1+
[X,\der_1]\cdot \der_2,$$
there are two possibilities:
\begin{itemize}
\item
if $1,2\in \{\sigma(1),\sigma(2),\sigma(3)\}$ or
$1,2\in \{\sigma(4),\sigma(5),\sigma(6)\}$ then
$(\sigma(4),\sigma(5),\sigma(6))=(1,2,s),$
and $(\sigma(1),\sigma(2),\sigma(3))=(q,r,6),$
where $\{q,r,s\}=\{3,4,5\},$ and $q<r.$
\item
if each of the following subsets
$\{\sigma(1),\sigma(2),\sigma(3)\}$
and $\{\sigma(4),\sigma(5),\sigma(6)\}$ contains exactly
one element $s\in \{1,2\}.$
\end{itemize}
Therefore,
$$
s_6(X_1,\ldots,X_6)=
s_6'(X_1,\ldots,X_6)+
s_6''(X_1,\ldots,X_6),$$
where
$$s_6'(X_1,\ldots,X_6)=$$
$$s_3^{rsym}(X_1,X_3,X_6)\circ s_3^{rsym}(X_2,X_4,X_5)-
s_3^{rsym}(X_1,X_4,X_6)\circ s_3^{rsym}(X_2,X_3,X_5)$$
$$+s_3^{rsym}(X_1,X_5,X_6)\circ s_3^{rsym}(X_2,X_3,X_4)
-s_3^{rsym}(X_2,X_3,X_6)\circ s_3^{rsym}(X_1,X_4,X_5)$$
$$+s_3^{rsym}(X_2,X_4,X_6)\circ s_3^{rsym}(X_1,X_3,X_5)-
s_3^{rsym}(X_2,X_5,X_6)\circ s_3^{rsym}(X_1,X_3,X_4)
,$$

$$s_6''(X_1,\ldots,X_6)=$$
$$-s_3^{rsym}(X_4,X_5,X_6)\circ s_3(X_1,X_2,X_3)+
s_3^{rsym}(X_3,X_5,X_6)\circ s_3(X_1,X_2,X_4)$$
$$-s_3^{rsym}(X_3,X_4,X_6)\circ s_3(X_1,X_2,X_5).
$$
Here we use notation $s_3^{rsym}$ instead of $s^{rsym.r}$ or $s_3^{rsym.l},$
because $s^{rsym.r}=s_3^{rsym.l}$ for any right-symmetric algebra.

Notice that
$$s_3(\der_1,\der_2,x_1\der_1)=\der_1\der_2,$$
$$s_3(\der_1,\der_2,x_2\der_1)=-\der_1^2,$$
$$s_3(\der_1,\der_2,x_1\der_2)=\der_2^2,$$
$$s_3(\der_1,\der_2,x_2\der_2)=-\der_1\der_2.$$
Therefore,
$$s_6''(X_1,\ldots,X_6)=$$
\begin{equation}\label{anuar}
s_3^{rsym}(x_1\der_1-x_2\der_2,x_1\der_2,X_6)\circ \der_1^2+
2 s_3^{rsym}(x_2\der_1,x_1\der_2,X_6)\circ \der_1\der_2
\end{equation}
$$-s_3^{rsym}(x_2\der_1,x_1\der_1-x_2\der_2,X_6)\circ \der_2^2
$$

Now calculate  $s_6''(X_1,\ldots,X_6)$ for $X_6=x_1^2\der_2.$
By (\ref{anuar}) we have
$$s_6''(X_1,\ldots,X_5,x_1^2\der_2)=$$
$$2 x_1^2\der_2\circ \der_1\der_2+(4x_1^2\der_1+6 x_1 x_2\der_2)\circ
\der_2^2=0.$$

Check that $s_6'(X_1,\ldots,X_6)=0$ for $X_6=x_1^2\der_2.$

Let $a,b,c\in <x_2\der_1, x_1\der_1-x_2\der_2,x_1\der_2\rangle.$ Notice
that
$$s_3^{rsym}(\der_i,a,X_6)=a\circ [X_6, \der_i]+X_6\circ[\der_i,\circ a]=$$
$$a\circ \der_i(X_6)-X_6\circ \der_i(a),$$

$$s_3^{rsym}(\der_j,b,c)=b\circ [c,\der_j]+c\circ [\der_j,b]=b\circ\der_j(c)
-c\circ \der_j(b).$$
By these formulas, it is easy to calculate that
$$
s_3^{rsym}(X_1,X_3,X_6)=s_3(\der_1,x_2\der_1,x_1^2\der_2)=
2 x_2\der_1\circ x_1\der_2=2x_1\der_1,$$

$$s_3^{rsym}(X_2,X_4,X_5)=s_3(\der_2,x_1\der_1-x_2\der_2,x_1\der_2)=
x_1\der_2\circ \der_2=0,$$
and
$$s_3^{rsym}(X_1,X_3,X_6)\circ s_3(X_2,X_4,X_5)=0.$$
Furthermore,
$$s_3^{rsym}(X_1,X_4,X_6)=s_3^{rsym}(\der_1,x_1\der_1-x_2\der_2,x_1^2\der_2)=$$
$$2 (x_1\der_1-x_2\der_2)\circ x_1\der_2-x_1^2\der_2\circ \der_1=$$
$$-2 x_1\der_2-2x_1\der_2=-4 x_1\der_2,$$

$$s_3^{rsym}(X_2,X_3,X_5)=s_3^{rsym}(\der_2,x_2\der_1,x_1\der_2)=
-x_1\der_2\circ \der_1=-\der_2,$$
and
$$s_3^{rsym}(X_1,X_4,X_6)\circ s_3(X_2,X_3,X_5)=0.$$
Finally,
$$s_3^{rsym}(X_1,X_5,X_6)=s_3(\der_1,x_1\der_2,x_1^2\der_2)=$$
$$2 x_1\der_2\circ x_1\der_2-x_1^2\der_2\circ \der_2=0,$$
and
$$s_3^{rsym}(X_1,X_5,X_6)\circ s_3^{rsym}(X_2,X_3,X_4)=0.$$
Similarly,
$$s_3^{rsym}(X_2,X_3,X_6)\circ s_3^{rsym}(X_1,X_4,X_5)=0,$$
$$s_3^{rsym}(X_2,X_4,X_6)\circ s_3^{rsym}(X_1,X_3,X_5)=0,$$
$$s_3^{rsym}(X_2,X_5,X_6)\circ s_3^{rsym}(X_1,X_3,X_4)=0.$$

So, we have established that $s_6'(X_1,\ldots,X_6)=0$ for $X_6=x_1^2\der_2.$
Thus,
$$s_6(X_1,\ldots,X_5,x_1^2\der_2)=
s_6'(X_1,\ldots,X_5,x_1^2\der_2)+
s_6''(X_1,\ldots,X_5,x_1^2\der_2)=0.$$

Similarly, one can calculate that
$$s_6(X_1,\ldots,X_5,X_6)=0,$$
for any $X_6=x_1^2\der_1-2x_1x_2\der_2, x_2^2\der_2-2x_1x_2\der_1,
x_2^2\der_1.$ In other words,
$esc(s_6)(X_1,\ldots,X_6)=0,$ for any
$(X_1,\ldots,X_6)\in V.$

Therefore, by (\ref{escort}) $s_6=0$
is an identity on $Vect_0(2).$

\section{$s_7=0$ is an identity on $Vect(2).$}

\begin{lm}\label{arman}
$s_7=0$ is an identity on $Vect(2).$
\end{lm}

{\bf Proof.} We see that $esc(s_7)$ is uniquely defined by the
homomorphism of $sl_2$-modules $f:L_1\rightarrow L_{-1}$ given by
$f(X)=s_7(\der_1,\der_2,x_1\der_1,x_2\der_2,x_1\der_1,x_2\der_2,X).$
Since
$$L_1\cong R(\pi_1)+R(2\pi_1+\pi_2),$$
and this isomorphism of $sl_2$-modules
can be given by divergence map,
it is clear that $f(X)=\lambda Div(X)$ for some $\lambda\in {\bf C}.$
Using the decomposition
$s_7^{rsym.r}=s_4^{rsym.r}\overset{\circ}\smile s_3,$
one can calculate that
$$s_7(\der_1,\der_2,x_1\der_1,x_2\der_2,x_1\der_1,x_2\der_2,x_1^2\der_1)=0.$$
Thus, $\lambda=0$ and $s_7=0$ is an identity on $Vect(2).$

\section{$5$- and $6$-commutators are primitive.}

The commutator operation is an invariant operation with two arguments
on the space of vector fields.
Therefore, any iteration of commutator operations $k-1$ times
gives us a new invariant operation with $k$ arguments on the space
of vector fields. Any linear combination of invariant $k$-operations
is also a $k$-operation. Below we show that
for $L=Vect_0(2)$ the $5$-commutator $s_5$ can not be obtained in a such way
from commutator operations. So, $5$-commutator $s_5$ is a new invariant
operation on $Vect_0(2)$ that can not be reduced to commutators.
Similar results hold for $6$-commutator.

\begin{lm} \label{19fev}
There does not exist a Lie polynomial $f=f(t_1,\ldots,t_5)$  such that
$s_5(X_1,\ldots,X_5)=f(X_1,\ldots,X_5),$ for any
$X_1,\ldots,X_5\in Vect_0(2).$
Similarly, one can not represent a $6$-commutator on $Vect(2)$ in the form
$s_6(X_1,\ldots,X_6)=g(X_1,\ldots,X_6),$ for any
$X_1,\ldots,X_6\in Vect(2),$ where $g$ a is Lie polynomial in $6$ variables.
\end{lm}

{\bf Proof.}
Let $L$ be a Lie algebra, $U(L)$  its universal enveloping algebra and
$$\Delta: U(L)\rightarrow U(L)\otimes U(L),\quad \Delta(X)=X\otimes 1+1\otimes X,\;\; \forall X\in L,$$
a comultiplication. For any $X_1,\ldots,X_k\in L,$
$$
\Delta(X_1\cdot \ldots \cdot X_k)=\sum_{l=0}^k
\sum_{\sigma\in Sym_{l,k-l}}X_{\sigma(1)}\cdots \ldots \cdot X_{\sigma(k-l)}
\otimes X_{\sigma(l+1)}\cdot \ldots\cdot X_{\sigma(k)}.$$
Thus, for any $X_1,\ldots,X_k\in L,$
$$
\Delta(s_k(X_1,\ldots,X_k))=\sum_{l=0}^k\sum_{\sigma\in Sym_{l,k-l}}
s_l(X_1,\ldots,X_l)\otimes s_{k-l}(X_{l+1},\ldots,X_{k}).$$
Therefore, if $s_k$ is the standard $k$-commutator,
i.e., if $s_k$ is obtained from Lie polynomial,
then  \cite{Jacobson}
$$
G_k=\sum_{l=1}^{k-1}s_l(X_1,\ldots,X_l)\otimes s_{k-l}(X_{l+1},\ldots,X_k).
$$
should be identically 0 for any $X_1,\ldots,X_k\in L.$ Here
$L=W_2$ if $k=6,$ and $L=S_1$ if $k=5.$

In a calculation of $G_k$ below we use  formulas for quadratic
parts of $k$ commutators (lemmas \ref{quadratic part for s3},
\ref{quadratic part for s4},
\ref{quadratic part for s5}).

Consider the case of $5$-commutators. Take $$(X_1,X_2,X_3,X_4,X_5)=
(\der_1,\der_2,x_1\der_1-x_2\der_2,x_2\der_1,x_1\der_2).$$
One can calculate that
$$G_5=$$
$$
-4\der_1\otimes \der_2
-4\der_2\otimes \der_1
-2\der_2\otimes x_1\der_1^2
-4\der_2\otimes x_2\der_1 \der_2
+4\der_1\der_2\otimes x_1\der_1
$$
$$
-4\der_1\der_2\otimes x_2\der_2
+4x_1\der_1\otimes \der_1\der_2
-2x_1\der_1^2\otimes \der_2
-4x_2\der_2\otimes \der_1\der_2
-4x_2\der_1\der_2\otimes\der_2
$$
$$\ne 0.$$
So, $s_5$ on $S_1$ can not be obtained from any Lie polynomial.

Consider now the case of $6$-commutator. Take
$$(X_1,X_2,X_3,X_4,X_5,X_6)=
(\der_1,\der_2,x_1\der_1,x_2\der_1,x_1\der_2,x_1^2\der_1).$$
We see that
$$s_3(X_1,X_2,X_3)\otimes s_3(X_4,X_5,X_6)=
\der_1\der_2\otimes(
3 x_1^2\der_1+2x_1x_2\der_2+x_1^3\der_1^2+
2 x_1^2x_2\der_1\der_2),$$
Therefore, $G_6$ has the term of the form
$\der_1\der_2\otimes x_1^3\der_1^2.$
Collect all terms of $G_6$ of the form
$\lambda \der_1\der_2\otimes x_1^3\der_1^2.$
Then their sum, denoted  by
$R,$ should be $0$ if $s_6$ is standard $5$-commutator.
As a differential operator of second order, $x_1^3\der_1^2$ can not
appear in  $s_2(X_i,X_j).$
Direct calculations then show that the elements  of the form
$s_l(X_{j_1},\ldots,X_{j_l}), j_1<\cdots<j_l, l=3,4,5,$ may
have the part $\mu x_1^3\der_1^2, \mu\ne 0$ only in one case:  $l=3,
(j_1,j_2,j_3) =(4,5,6).$ So, $R=
\der_1\der_2\otimes x_1^3\der_1^2\ne 0.$
This contradiction shows that $6$-commutator on $W_2$
is primitive.

\section{$s_5$ and $s_6$ are cocycles}

Let $d: C^k(L,L)\rightarrow C^{k+1}(L,L)$ be the coboundary operator.
Then $$d\psi=d'\psi+d''\psi,$$ where
$$d'\psi(X_1,\ldots,X_{k+1})=\sum_{i<j}(-1)^{i+j}\psi([X_i,X_j],X_1,\ldots,\hat{X_i},\ldots,\hat{X_j},\ldots,X_{k+1}),$$
$$d''\psi(X_1,\ldots,X_{k+1})=\sum_{i=1}^{k+1}(-1)^{i+1}[X_i,\psi(X_1,\ldots,\hat{X_i},\ldots,X_{k+1}).$$

\begin{lm}\label{fev25}
$d's_k^{rsym.r}=0,$ if $k$ is even and $d's_k^{rsym.r}=-s_{k+1}^{rsym.r},$
if $k$ is odd.
\end{lm}

{\bf Proof.} This follows from induction in $n$ and the following
relation
$$s_{k+1}^{rsym.r}=\sum_{i=1}^{k+1}(-1)^{i+k+1}
(s_k^{rsym.r}(X_1,\ldots,\hat{X_i},\ldots,X_{k+1}))X_i.$$

\begin{lm}\label{zamok versal} $(2d'+d'')s_{k}=0,$
for any $k\ge n^2+2n-2.$
\end{lm}

{\bf Proof.} By corollary~\ref{alia}, $ad\,X\in Der(Vect(n),s_k),$
if $k\ge n^2+2n-2.$ Therefore,
$$[X_i,s_{k}(X_1,\ldots,\hat{X_i},\ldots,X_{k+1})]=$$
$$\sum_{j=1}^{i-1}(-1)^{1+j}s_k([X_i,X_j],\ldots,\hat{X_i},\ldots,\hat{X_j},\ldots,X_{k+1})$$
$$+\sum_{j=i+1}^{k+1}(-1)^{j}s_k([X_i,X_j],\ldots,\hat{X_i},\ldots,\hat{X_j},\ldots,X_{k+1}),$$
and
$$d''s_k(X_1,\ldots,X_{k+1})=$$
$$-2\sum_{i<j}(-1)^{i+j}s_k([X_i,X_j],\ldots\hat{X_i},\ldots,\hat{X_j},\ldots,X_{k+1}).$$
In other words, $d''s_k=-2 d's_k,$ if $k\ge n^2+2n-2.$

\begin{crl} \label{04march} $ds_5=0$ on $Vect_0(2).$
\end{crl}

{\bf Proof.} By lemma \ref{bir} $s_5=s_5^{rsym.r}$ and $s_6=s_6^{rsym.r}=0$
are identities on $Vect_0(2).$ Therefore, by lemma
$d_{s_5}=d's_5+d''s_5=-d's_5=s_6=0$ is an identity on $Vect_0(2).$

\begin{crl}\label{march04} $ds_6=0$ on $Vect(2).$
\end{crl}

{\bf Proof.}
By lemma \ref{28fevruary} and lemma \ref{arman} $s_6=s_6^{rsym.r}$ and
$s_7=s_7^{rsym.r}$ is an identity on $Vect(2).$
Therefore, by lemma \ref{fev25} and \ref{zamok versal}
$d_{s_6}=d's_6+d''s_6=-d's_6=0$ is an identity on $Vect(2).$

{\bf Remark.} One can prove that
$(L,\{s_2,s_{l}\})$ is also $sh$-Lie, for $l=n^2+2n-2,$ if $L=Vect(n)$ and
$l=n^2+2n-3,$ if $L=Vect_0(n).$

Our results can be formulated in terms of generalized cohomology operators.
There are two ways to do it. In the first way
one saves the index of nilpotency  $d^2=0,$ but changes the grading degree.
In the second way one saves grading degree, but changes the index of
nilpotency from $d^2=0$ to $d^N=0.$ A cohomology theory for $d^N=0$ was
developed in \cite{Dubois-Violette}.

Let us show how to do it for left multiplication operators. Let
$L=Vect(n)$ be the right-symmetric algebra of vector fields and
$l_a$ denotes left multiplication operator, $(b)l_a= a\circ b.$
Define a linear operator $d:\wedge ^*(L,L)\rightarrow \wedge
^*(L,L)$ by
$$d:C^k(L,L)\rightarrow C^{k+n}(L,L),$$
$$d\psi(a_1,\ldots,a_{k+n})=\sum_{\sigma\in Sym_{n,k}}
\,sign\,\sigma\, l_{a_{\sigma(1)}}\cdots l_{a_{\sigma(n)}}
\psi(a_{\sigma(n+1)},\ldots,a_{\sigma(k+n)}).$$
Then the condition $d^2=0$ follows from theorem 3.3 of \cite{Dzh00}.

In the second case we need to consider a coboundary operator with
grading degree $+1,$
$$d_l: \wedge^*(L,L)\rightarrow \wedge^*(L,L),$$
$$d_l: \wedge^{k}(L,L)\rightarrow \wedge^{k+1}(L,L),$$
$$d_l\psi(a_1,\ldots,a_{k+1})=\sum_{i=1}^{k+1}
(-1)^il_{a_i}\psi(a_1,\ldots,\hat{a_i},\ldots,a_{k+1}).$$
Then $d_l^{2n}=0.$

One can construct similar coboundary operators  corresponding to
right-multiplication operators. For example,
$$d_r: \wedge^*(L,L)\rightarrow \wedge^*(L,L),$$
$$d_r: \wedge^{k}(L,L)\rightarrow \wedge^{k+1}(L,L),$$
$$d_r\psi(a_1,\ldots,a_{k+1})=\sum_{i=1}^{k+1}
(-1)^ir_{a_i}\psi(a_1,\ldots,\hat{a_i},\ldots,a_{k+1}).$$
has the property $d_r^{n^2+2n-1}=0.$

These constructions have some other modifications that include
the case of more general right-symmetric algebras and their modules.

\section{Proofs  of main results}

{\bf Proof of theorem \ref{N-commutator}.} This follows from lemmas
\ref{28fevruary}, \ref{fev28}, \ref{n^2+2n-3 left commutative} and
corollary \ref{n^2+2n commutator}.

{\bf Proof of theorem \ref{5-commutator}.} This follows from lemmas
\ref{4-left commutative}, \ref{s5 for adjoint}, \ref{bir},
\ref{19fev} and corollary \ref{04march}.

{\bf Proof of theorem \ref{6-commutator}.} This follows from lemmas
\ref{arman}, \ref{19fev},\ref{5-left commutative} and corollary
\ref{march04}.

\begin{center}
{\em ACKNOWLEDGEMENTS}
\end{center}

I am grateful to INTAS foundation for support. I wish to  thank
M.~Kontsevich and J.~Hoppe for the stimulating interest these
results.

\end{document}